\documentclass[11pt]{article}
\usepackage{amsmath}
\usepackage{amssymb}
\usepackage{mathrsfs}
\usepackage{amsthm}
\usepackage{microtype}
\usepackage{needspace}
\usepackage[colorlinks=true,citecolor=blue,linkcolor=blue,urlcolor=blue]{hyperref}
\hypersetup{
 pdftitle={Spectral Schwarz Lemmas and Parallel Obstructions to Metric Comparison},
 pdfauthor={Xu Pang, Yao Sha, Ren Xin-An},
 pdfsubject={Complex differential geometry; manuscript version V13},
 pdfkeywords={Schwarz lemma, holomorphic bundle map, singular value, Schatten moment, rank rigidity}
}
\theoremstyle{plain}
\newtheorem{theorem}{Theorem}[section]
\newtheorem{lemma}[theorem]{Lemma}
\newtheorem{prop}[theorem]{Proposition}
\newtheorem{corollary}[theorem]{Corollary}
\theoremstyle{definition}
\newtheorem{defi}[theorem]{Definition}
\newtheorem{example}[theorem]{Example}
\theoremstyle{remark}
\newtheorem{remark}[theorem]{Remark}
\newtheorem*{quotedtheorem}{Theorem}

\begin{document}
	\title{\bf Spectral Schwarz Lemmas and Parallel Obstructions to Metric Comparison}
	\author{Xu Pang, Yao Sha, Ren Xin-An}

	
	\date{}
	
	{\setlength{\tabcolsep}{0pt}\maketitle}
	
	\renewcommand{\theequation}{\thesection.\arabic{equation}}
	
	\numberwithin{equation}{section}

\begin{abstract}
In this paper, we characterize when comparison of contracted Chern curvatures implies comparison of Hermitian bundle metrics with the same constant on compact Hermitian bases.
If the curvature contraction of the reference metric is semipositive, this implication holds for every Hermitian metric
if and only if the kernel of this contraction contains no nonzero subbundle parallel with respect to the Chern connection of the reference metric.
The equality terms of a $C^2$ spectral truncation identify the
parallel summand on which comparison can fail, without regularity assumptions on the curvature kernel.
We also establish a sharp operator estimate, with a constant independent of the
bundle ranks, for holomorphic bundle maps on complete K\"ahler bases under
curvature bounds and bounded base-map energy density.
Smooth Schatten moments give global control without a separate
boundedness assumption on the bundle map.
Finally, we obtain rank bounds and conditions for vanishing on compact
Hermitian bases from directional curvature traces over fiber
subspaces.
For higher-dimensional subspaces, these conditions allow mixed
signs on individual fiber lines.
\end{abstract}



\section{Introduction}

In its Schwarz--Pick form, the classical Schwarz lemma states that holomorphic self-maps of the unit disk decrease the Poincar\'{e} metric.
Ahlfors extended this principle to negatively curved metrics on Riemann surfaces \cite{Ahl38}, and Chern developed higher-dimensional versions for Hermitian manifolds \cite{Che68}.
Yau proved a Schwarz lemma for holomorphic maps from complete K\"{a}hler manifolds by combining a Bochner formula with a generalized maximum principle \cite{Yau78}.
His estimate bounds the energy density in terms of a lower Ricci bound on the domain and a negative upper bisectional curvature bound on the target.

Royden weakened Yau's target assumption to a holomorphic sectional curvature bound for K\"ahler targets \cite[Theorem~1]{Roy80}. Chen--Cheng--Lu \cite{CCL79} and Royden
\cite[Theorem~2]{Roy80} obtained pointwise upper bounds for the pullback of the target metric in terms of the domain metric, with constants determined by the holomorphic sectional curvature bounds on both manifolds. On complete domains, these estimates require additional lower curvature bounds, which are automatic in the compact case.
In particular, if $(M,g)$ is compact K\"ahler, $(N,h)$ is K\"ahler, and their holomorphic sectional curvatures satisfy
$\operatorname{HSC}(M,g)\ge-K$ and $\operatorname{HSC}(N,h)\le-B<0$, where $K\ge0$, then every holomorphic map $f:M\to N$ satisfies
\[
 f^*\omega_h\le\frac{K}{B}\omega_g.
\]
Ni gave a spectral treatment of this comparison \cite[Theorem~1.4]{Ni21}. Beyond the K\"ahler setting, related work develops Bochner formulas and curvature conditions for Hermitian Schwarz lemmas \cite{Tang21,BS26,YZ19} and extensions to almost-Hermitian manifolds \cite{Tos07}. For broader background, see \cite{Bro22,Kob98}.

Yau's energy estimate bounds the sum of the squared singular values of $df$, whereas pointwise metric comparison bounds the largest one.
Different curvature contractions and eigenvalue weights arise in the associated Bochner formulas, so the choice of spectral quantity influences which curvature bounds can be used.

For two Hermitian metrics $G,H$ on a holomorphic bundle $E\to X$, the squared singular values of the identity map $(E,G)\to(E,H)$ are the eigenvalues of $G^{-1}H$. Bounding the largest eigenvalue by $\kappa$ is therefore equivalent to the metric comparison $H\le\kappa G$. Fix a Hermitian base metric $b$ and write
\[
 S_b(G)=\operatorname{tr}_{\omega_b}R^G,
 \qquad
 S_b(H)=\operatorname{tr}_{\omega_b}R^H.
\]
Each curvature tensor has its fiber indices lowered by its own metric, so these contractions are Hermitian forms on $E$.
The question is when comparison of these forms implies comparison of the metrics.
The following theorem of Wang--Yang--Yau gives such a conclusion when $S_b(G)>0$.

\begin{quotedtheorem}[{Wang--Yang--Yau
\cite[Theorem~1.2 and Corollary~3.4]{WYY26}}]
Let $(X,b)$ be compact K\"ahler, and let $G,H$ be Hermitian metrics on a holomorphic vector bundle $E\to X$. If
\[
 S_b(G)>0,\qquad S_b(H)\le\kappa S_b(G),\qquad \kappa>0,
\]
then $H\le\kappa G$.
\end{quotedtheorem}

Xiong--Yang--Yau's directional comparison \cite[Theorem~1.8]{XYY} also implies this conclusion on a compact Hermitian base. The positive contracted bound supplies a common tangent direction for the two curvature tests, as explained in Remark~\ref{rem:metric-xyy}.

Under semipositivity, a $G$-parallel holomorphic summand annihilated by $S_b(G)$ obstructs metric comparison:
a constant enlargement of the metric on this summand can preserve the curvature inequality while violating $H\le\kappa G$.
The pointwise kernel of $S_b(G)$, however, need not have constant rank or be preserved by the Chern connection.
Thus this construction leaves open whether comparison can fail without such a parallel summand in the kernel.
Our first theorem shows that every failure has precisely this structure: the eigenspaces of $G^{-1}H$ above $\kappa$
form a parallel holomorphic summand contained in the curvature kernel, and the relative endomorphism is parallel there.

\Needspace{15\baselineskip}
\begin{theorem}
\label{thm:excess-main}
Let $(X,b)$ be compact connected Hermitian, and let $G,H$ be Hermitian metrics on a holomorphic vector bundle $E\to X$. Suppose
that
\[
 S_b(G)\ge0,\qquad S_b(H)\le\kappa S_b(G),\qquad \kappa>0.
\]
Let $A=G^{-1}H$, defined by $H(\xi,\bar\eta)=G(A\xi,\bar\eta)$, and let $(E_+)_x$ be the sum of the eigenspaces of $A_x$ with eigenvalues greater than $\kappa$.
Then:
\begin{enumerate}
\item $E_+$ is a smooth holomorphic subbundle, and
$E=E_+\oplus E_+^{\perp_G}$ is preserved by the Chern connection of $G$.
\item $S_b(G)$ annihilates $E_+$, and $A|_{E_+}$ is parallel.
In particular, all eigenvalues above $\kappa$ are constant.
\end{enumerate}
Consequently, $H\le\kappa G$ if $E$ has no nonzero $G$-parallel subbundle annihilated by $S_b(G)$.
\end{theorem}

In particular, the assumption $S_b(G)>0$ in Wang--Yang--Yau's
comparison theorem can be weakened to semipositivity everywhere
and positive definiteness at a single point
(Corollary~\ref{thm:metric-comparison}).
More generally, for fixed $G$ and $\kappa$, comparison holds for
\emph{every} metric $H$ satisfying
$S_b(H)\le\kappa S_b(G)$ if and only if
$\ker S_b(G)$ contains no nonzero $G$-parallel subbundle
(Corollary~\ref{cor:universal-comparison}).
This criterion can hold even when $S_b(G)$ is degenerate everywhere.

A concrete example is the pullback of the universal quotient bundle on $\mathbb{CP}^2$ by the quadratic Veronese embedding of $\mathbb{CP}^1$. The bundle has rank two, and its quotient metric $G$ has semipositive contracted curvature of rank one at every point. Its contracted curvature kernel is a smooth line subbundle, but is not preserved by the Chern connection of $G$. Hence $S_b(H)\le\kappa S_b(G)$ implies $H\le\kappa G$ (Example~\ref{ex:veronese-comparison}). The strict positivity assumptions of the two comparisons cited above fail for this metric. On compact K\"ahler manifolds, the theorem gives metric comparison under $\operatorname{Ric}_g\ge0$ when there is no nonzero Ricci-flat parallel tangent summand; local irreducibility and $\operatorname{Ric}_g\not\equiv0$ suffice (Corollary~\ref{cor:semipositive-tangent-comparison}).

Wang--Yang--Yau's proof uses spectral weights concentrating near the largest relative eigenvalue \cite[Proof of Theorem~3.1]{WYY26}.
Their contradiction argument relies on strict positivity of the contraction $S_b(G)$.
Here we apply the spectral Bochner formula to
\[
 F=\operatorname{tr}(A-\kappa I)_+^3,
\]
where the positive part acts on the eigenvalues.
The cubic truncation is $C^2$ at the threshold, and its trace remains $C^2$ at eigenvalue crossings.
Under semipositivity, the curvature inequality gives $\Delta_bF\ge0$, so $F$ is constant on the compact base.
The equality terms then force $(A-\kappa I)_+^3$ to be parallel. Its image gives the required splitting.
This recovers the entire excess spectral subbundle without assuming a spectral gap or regularity of the curvature kernel.

We next consider global operator bounds for general holomorphic bundle maps on complete K\"ahler bases.
The largest squared singular value need not attain a maximum or be bounded a priori, so global control must be established from the curvature inequalities.
Let $f:(M,g)\to(N,h)$ be holomorphic, and let $(E_1,G)\to M$ and $(E_2,H)\to N$ be Hermitian holomorphic vector bundles of ranks $r_1,r_2$. A holomorphic bundle map
$\varphi:E_1\to E_2$ covering $f$ is equivalently a holomorphic section of $E_1^*\otimes f^*E_2$.
Write $(\varphi^*H)(\xi,\bar\eta)=H(\varphi\xi,\overline{\varphi\eta})$ for the induced Hermitian form. Our next theorem bounds $\varphi^*H$ by a multiple of $G$,
with a constant independent of the bundle ranks.

\Needspace{13\baselineskip}
\begin{theorem}
\label{thm:complete-main}
Let $(M,g)$ be complete K\"ahler with Ricci curvature bounded from below, and let $(N,h)$ be Hermitian. Let $f:M\to N$ be holomorphic and let $\varphi:(E_1,G)\to(E_2,H)$ be a holomorphic bundle map covering $f$. Assume
\[
 E_f:=\sup_M\operatorname{tr}_{\omega_g}f^*\omega_h<\infty.
\]
If $K_1,K_2>0$ satisfy
\[
 \operatorname{tr}_{\omega_g}R^{E_1}\ge K_1\varphi^*H,
 \qquad
 R^{E_2}(w,\bar w,\eta,\bar\eta)
 \le K_2|w|_h^2|\eta|_H^2
\]
for all $w\in T^{1,0}N$ and $\eta\in E_2$ over the same point, then
\[
 \varphi^*H\le\frac{K_2E_f}{K_1}G.
\]
\end{theorem}

The estimate is sharp on both compact and complete noncompact bases (Remark~\ref{rem:sharp-complete}).
Under the same curvature bounds, Xiong--Yang--Yau obtained, for nonconstant $f$ and nonzero $\varphi$ on compact Hermitian bases, the trace estimate
\[
 \sup_M\operatorname{tr}_G\varphi^*H
 \le r_\varphi\frac{K_2E_f}{K_1},
 \qquad
 r_\varphi=\max_M\operatorname{rank}_{\mathbb C}\varphi
\]
\cite[Theorem~1.13]{XYY}.
Our operator estimate recovers this trace bound by summing over the nonzero eigenvalues, while avoiding the factor $r_\varphi$ incurred when the trace is used to bound the largest eigenvalue. On compact Hermitian bases, the same operator estimate follows from the local-support proof of Proposition~\ref{thm:traced-sum} with $k=1$.

To obtain global control on a complete base without a separate boundedness assumption on $\varphi$, we use the smooth moments
\[
 A=G^{-1}\varphi^*H,\qquad U_p=\operatorname{tr}(A^p).
\]
The lower bound on the curvature contraction of $E_1$ contributes a positive multiple of $\operatorname{tr}(A^{p+1})$ in the Bochner inequality for $U_p$.
H\"older's inequality and a nonlinear Omori--Yau principle give, for every integer $p\ge2$,
\[
 \sup_M U_p\le q\left(\frac{K_2E_f}{K_1}\right)^p,
 \qquad q=\min\{r_1,r_2\}.
\]
Taking $p$-th roots and letting $p\to\infty$ yields the rank-independent operator estimate, since $q^{1/p}\to1$.

On a compact base, we also study distortion in several fiber directions through partial sums and products. Let
$|\lambda_1|\ge\cdots\ge|\lambda_q|$ denote the singular values of $\varphi$ with respect to $G$ and $H$. Write $\mu_a=|\lambda_a|^2$ for the largest $q$ eigenvalues of $A$. Set
\[
 \sigma_k=\sum_{a=1}^k\mu_a,
 \qquad
 \pi_k=\prod_{a=1}^k\mu_a
       =\|\wedge^k\varphi\|_{\mathrm{op}}^2.
\]
We adapt Ni's local-support method \cite{Ni19} to general holomorphic bundle maps. On compact Hermitian bases, we obtain partial-sum estimates directly from these supports, recovering the operator bound when $k=1$. Our product estimates require only curvature traces over $k$-dimensional fiber subspaces, rather than bounds on individual fiber directions.
These weaker hypotheses allow mixed signs in the fibers. For example, on $\mathcal O(-1)\oplus\mathcal O(3)$ over $\mathbb{CP}^1$, they give $\pi_2\le1$ although the largest
squared singular value is $3/2$ (Example~\ref{ex:product-mixed-sign}).
Thus the product estimate need not follow from metric comparison with the same constant.

Curvature traces over $l$-planes can force $\wedge^l\varphi$ to vanish even when the curvature on individual fiber lines has no fixed sign.
For $l>1$, this restricts the rank without forcing $\varphi$ itself to vanish. For a complex $l$-plane $\sigma\subset(E_1)_x$, define
\[
 R^{E_1}(v;\sigma)
 =\sum_{a=1}^l R^{E_1}(v,\bar v,e_a,\bar e_a),
\]
where $e_1,\ldots,e_l$ is any $G$-unitary basis of $\sigma$. Define the target trace in an $H$-unitary basis.

\Needspace{12\baselineskip}
\begin{theorem}\label{thm:intro-rank}
Let $(M,g)$ and $(N,h)$ be Hermitian, with $M$ compact. Let $f:M\to N$ be holomorphic and let $\varphi:(E_1,G)\to(E_2,H)$ be a holomorphic bundle map covering $f$. Fix $1\le l\le q$. Suppose that for every $x\in M$ and every pair of complex $l$-planes $\sigma\subset(E_1)_x$, $\sigma'\subset(E_2)_{f(x)}$, some nonzero $v\in T_x^{1,0}M$ satisfies
\[
 R^{E_1}(v;\sigma)\ge0,\qquad
 R^{E_2}(df_x(v);\sigma')\le0,
\]
with at least one strict inequality. Then $\wedge^l\varphi=0$, or equivalently,
$\operatorname{rank}_{\mathbb C}\varphi_x<l$ at every point.
\end{theorem}

The same base direction must satisfy both curvature tests, but it may depend on the two fiber planes.
The hypotheses are imposed on fiber $l$-planes, equivalently on decomposable exterior vectors, without a sign assumption
on each fiber line. A direct application of the linewise vanishing criteria of
Xiong--Yang--Yau to $\wedge^l\varphi$ would require the
corresponding hypotheses on all nonzero exterior vectors
\cite{XYY}.
When the direction depends on the planes, this extension does not follow merely by linearity.

The proof uses Ni's logarithmic spectral supports \cite{Ni19,Ni21}.
Dimension conditions produce the required common direction (Corollary~\ref{cor:dimension}). Under the corresponding non-strict trace signs in every base
direction, and assuming $r_1\le r_2$, the top spectral product $\pi_{r_1}$ is constant. In particular, full rank at one point implies full rank everywhere
(Theorem~\ref{thm:det-rigidity}).

The paper is organized as follows. Section~\ref{sec:preliminaries} establishes the global spectral formulas.
Sections~\ref{sec:semipositive} and~\ref{sec:complete} prove Theorems~\ref{thm:excess-main} and~\ref{thm:complete-main}.
Section~\ref{sec:compact} develops the local supports and compact
spectral estimates, and Section~\ref{sec:rank} treats determinant
rigidity and proves Theorem~\ref{thm:intro-rank}.
\section{Spectral Bochner formulas}\label{sec:preliminaries}

All manifolds are connected and without boundary. Dimensions and ranks are complex unless stated otherwise. We use the Chern
curvature convention for which, in a normal holomorphic frame at a point,
\[
 R^E_{i\bar j a\bar b}=-\partial_i\partial_{\bar j}H_{a\bar b}.
\]
With this convention, the Chern curvature of $(T^{1,0}\mathbb{CP}^n,g_{\mathrm{FS}})$ has positive holomorphic bisectional curvature.
For a Hermitian metric $g$ we write $\Delta_g u=g^{i\bar j}\partial_i\partial_{\bar j}u$ for the complex Laplacian. For real-valued $u$,
let $|\partial u|_g^2=g^{i\bar j}(\partial_i u)(\partial_{\bar j}u)$.
Thus $\Delta_g F(u)=F'(u)\Delta_g u+F''(u)|\partial u|_g^2$.
These conventions fix the signs and normalizing factors below. Hermitian forms are linear in the first argument.
The notation $A=G^{-1}B$ for a relative endomorphism is intrinsic:
$B(\xi,\bar\eta)=G(A\xi,\bar\eta)$. In matrix computations we distinguish this endomorphism from products of the ordinary Gram matrices of the Hermitian forms.

We first establish the global moment and spectral differentiation formulas.
The local supports for partial sums and products, with their Bochner formulas, are constructed in Section~\ref{sec:compact}.

\subsection{Conventions and spectral quantities}

Let $(M,g)$ and $(N,h)$ be Hermitian manifolds of dimensions $m$ and $n$, respectively. Let $(E_1,G)\to M$ and $(E_2,H)\to N$ be Hermitian holomorphic vector bundles of ranks $r_1$ and $r_2$. Let $f:M\to N$ be holomorphic and let $\varphi:E_1\to E_2$ be a holomorphic bundle map covering $f$.
Equivalently, $\varphi$ is a holomorphic section of
$E_1^*\otimes f^*E_2$. Set
\[
 q=\min\{r_1,r_2\},\qquad
 r_f=\max_M\operatorname{rank}_{\mathbb C}df,\qquad
 r_\varphi=\max_M\operatorname{rank}_{\mathbb C}\varphi.
\]
Fix $x_0\in M$ and choose holomorphic coordinates whose coordinate frames are unitary at $x_0$ and $f(x_0)$. Choose normal holomorphic
frames $\{e_a\}_{a=1}^{r_1}$ for $E_1$ near $x_0$ and $\{s_\gamma\}_{\gamma=1}^{r_2}$ for $E_2$ near $f(x_0)$, and pull back the latter to a neighborhood $B_{x_0}\subset M$.
The frames may be chosen so that
\begin{equation}\label{BJ}
\begin{aligned}
&\varphi(e_a)(x_0)=\lambda_a s_a\quad(1\le a\le q),\qquad
  \varphi(e_a)(x_0)=0\quad(q<a\le r_1),\\
&|\lambda_1|\ge\cdots\ge|\lambda_q|,\qquad
  \nabla e_a(x_0)=0,\qquad \nabla s_\gamma(f(x_0))=0,\\
&G(e_a,\bar e_b)(x_0)=\delta_{ab},\qquad
  H(s_\gamma,\bar s_\delta)(f(x_0))=\delta_{\gamma\delta}.
\end{aligned}
\end{equation}
Here and below $\nabla$ denotes the induced Chern connection. Write $\mu_a=|\lambda_a|^2$ for $1\le a\le q$ and $\mu_a=0$ for $q<a\le r_1$. These are the eigenvalues of
$\varphi^*H$ relative to $G$.
\begin{defi}\label{def:spectral}
For $1\le k\le q$, define
\[
 \sigma_k=\sum_{a=1}^k\mu_a,\qquad
 \pi_k=\prod_{a=1}^k\mu_a.
\]
\end{defi}
Both functions are continuous but may fail to be differentiable when eigenvalues meet. We have
\[
 \sigma_1=\mu_1,\qquad
 \sigma_q=\operatorname{tr}_G\varphi^*H,\qquad
 \varphi^*H\le\mu_1G\le\sigma_kG.
\]
Consequently any uniform upper bound for $\sigma_k$ gives a tensorial Schwarz estimate.
\begin{remark}
The induced Hermitian metrics on exterior powers give
\[
 \pi_k=\|\wedge^k\varphi\|_{\mathrm{op}}^2.
\]
Thus $\pi_k$ is the squared complex $k$-dilation. It measures the metric distortion on exterior powers; see \cite{Ni19,Ni21-2,Whi57}.
\end{remark}

We write $\varphi^\gamma_{ai}=(\nabla_i\varphi)_a^\gamma$; repeated coordinate indices are summed.
\subsection{The Schatten-moment formula}

We use $U_p$ for the globally defined Schatten moment
\begin{equation}\label{U}
 U_p=\operatorname{tr}(A^p)
     =\sum_{a=1}^{q}\mu_a^p,
 \qquad A=G^{-1}\varphi^*H,
\end{equation}
where $p$ is a positive integer. Yang--Chen \cite[p.~947, equation~(2.2)]{YC81} used these moments for a pullback tangent metric. We record the full formula for a general holomorphic bundle map, including the derivative terms.

\begin{prop}\label{prop:moment-bochner}
Let $p\ge1$ be an integer, let $U_p$ be defined by \eqref{U}, and fix $x_0\in M$.  In the normal frames \eqref{BJ}, with $\mu_a$ as above,
one has at $x_0$
\begin{align}
 \partial\bar\partial U_p
={}&p\sum_{a=1}^{r_1}\sum_{\gamma=1}^{r_2}
 \mu_a^{p-1}\varphi^\gamma_{ai}
 \overline{\varphi^\gamma_{aj}}\,
 dz^i\wedge d\bar z^j \notag\\
 &+p\sum_{a=1}^{r_1}\sum_{b=1}^{q}\sum_{\ell=0}^{p-2}
 \mu_a^\ell \mu_b^{p-1-\ell}
 \varphi^b_{ai}\overline{\varphi^b_{aj}}\,
 dz^i\wedge d\bar z^j \notag\\
 &+p\sum_{a=1}^{q}\mu_a^p
 \left(R^{E_1}_{i\bar ja\bar a}
 -f_i^\alpha\overline{f_j^\beta}
  R^{E_2}_{\alpha\bar\beta a\bar a}\right)
 dz^i\wedge d\bar z^j. \label{eq:moment-bochner}
\end{align}
A factor with exponent zero is understood as $1$, including when $\mu_a=0$. When $p=1$, the threefold sum is empty.
In particular, for every $v\in T_{x_0}^{1,0}M$,
\begin{equation}\label{eq:moment-ineq}
\partial\bar\partial U_p(v,\bar v)\ge
 p\sum_{a=1}^{q}\mu_a^p
 \left\{R^{E_1}(v,\bar v,e_a,\bar e_a)
 -R^{E_2}(f_*v,\overline{f_*v},s_a,\bar s_a)\right\}.
\end{equation}
\end{prop}

\begin{proof}
Set $B=\varphi^*H$. Let $\mathbf G=(G_{a\bar b})$ and $\mathbf B=(B_{a\bar b})$ be the ordinary Gram matrices, and put
$K=\mathbf G^{-1}\mathbf B$. This matrix has the same eigenvalues as the intrinsic relative endomorphism $A$: both spectra are determined
by $\det(\mathbf B-t\mathbf G)=0$. Thus
$U_p=\operatorname{tr}(K^p)$. Cyclic invariance of the trace gives
\[
 \partial_{\bar j}U_p
 =p\operatorname{tr}(K^{p-1}\partial_{\bar j}K).
\]
A second differentiation yields the exact identity
\begin{equation}\label{eq:matrix-power-hessian}
 \partial_i\partial_{\bar j}U_p
 =p\operatorname{tr}(K^{p-1}\partial_i\partial_{\bar j}K)
 +p\sum_{\ell=0}^{p-2}
  \operatorname{tr}\!\left(
  K^\ell(\partial_iK)K^{p-2-\ell}
  (\partial_{\bar j}K)\right).
\end{equation}

We first differentiate the pullback Hermitian form. In local frames,
\[
 B_{a\bar b}
 =H_{\gamma\bar\delta}(f)
   \varphi_a^\gamma\overline{\varphi_b^\delta}.
\]
For $1\le a\le q$, holomorphicity of $\varphi$, normality of the two frames, and the
curvature convention give
\begin{equation}\label{eq:B-second-derivative}
 \partial_i\partial_{\bar j}B_{a\bar a}
 =\sum_{\gamma=1}^{r_2}
   \varphi^\gamma_{ai}\overline{\varphi^\gamma_{aj}}
  -\mu_af_i^\alpha\overline{f_j^\beta}
   R^{E_2}_{\alpha\bar\beta a\bar a}.
\end{equation}
For $a>q$, the same computation has only its derivative sum, since $\varphi(e_a)(x_0)=0$.
At the singular frame, the first derivatives of $B$ are
\begin{equation}\label{eq:B-first-derivatives}
 \partial_iB_{a\bar b}
 =\varphi^b_{ai}\overline{\lambda_b},
 \qquad
 \partial_{\bar j}B_{b\bar a}
 =\lambda_b\overline{\varphi^b_{aj}}
 \quad(1\le a\le r_1,\ 1\le b\le q).
\end{equation}
For $b>q$, both displayed derivatives vanish.

We evaluate \eqref{eq:matrix-power-hessian} at $x_0$, where
$K=\operatorname{diag}(\mu_1,\ldots,\mu_{r_1})$ and normality of the
frame for $E_1$ gives
\begin{equation}\label{eq:A-first-derivatives}
 \partial_iK=\partial_i\mathbf B,
 \qquad
 \partial_{\bar j}K=\partial_{\bar j}\mathbf B.
\end{equation}
Using \eqref{eq:B-first-derivatives}, the second term on the
right-hand side of \eqref{eq:matrix-power-hessian} becomes
\begin{equation}\label{eq:moment-cross-term}
 p\sum_{a=1}^{r_1}\sum_{b=1}^{q}\sum_{\ell=0}^{p-2}
 \mu_a^\ell \mu_b^{p-1-\ell}
 \varphi^b_{ai}\overline{\varphi^b_{aj}}.
\end{equation}
Here, as dictated by $K^0=I$, a factor with exponent zero is understood to be $1$, also when the corresponding $\mu_a$ vanishes.

It remains to compute the diagonal entries of the second derivative of $K$.  Since the first derivatives of $\mathbf G$ vanish at $x_0$, differentiating
$\mathbf G^{-1}\mathbf G=I$ twice and using \eqref{eq:B-second-derivative} gives
\begin{equation}\label{eq:A-second-derivative}
 (\partial_i\partial_{\bar j}K)_{a\bar a}
 =\sum_{\gamma=1}^{r_2}
 \varphi^\gamma_{ai}\overline{\varphi^\gamma_{aj}}
 +\begin{cases}
 \mu_a\left(R^{E_1}_{i\bar ja\bar a}
 -f_i^\alpha\overline{f_j^\beta}
 R^{E_2}_{\alpha\bar\beta a\bar a}\right),&a\le q,\\
 0,&a>q.
 \end{cases}
\end{equation}
Thus the first term on the right-hand side of
\eqref{eq:matrix-power-hessian} equals
\begin{align*}
 {}&p\sum_{a=1}^{r_1}\sum_{\gamma=1}^{r_2}
 \mu_a^{p-1}\varphi^\gamma_{ai}
 \overline{\varphi^\gamma_{aj}}
 +p\sum_{a=1}^{q}\mu_a^p
 \left(R^{E_1}_{i\bar ja\bar a}
 -f_i^\alpha\overline{f_j^\beta}
 R^{E_2}_{\alpha\bar\beta a\bar a}\right).
\end{align*}
Combining this expression with \eqref{eq:moment-cross-term} proves \eqref{eq:moment-bochner}.  On evaluation at $(v,\bar v)$, each term in
the first two lines of \eqref{eq:moment-bochner} is a nonnegative coefficient times an absolute square.  Deleting those terms proves \eqref{eq:moment-ineq}.
\end{proof}

\begin{remark}\label{rem:moment-checks}
For $p=1$, formula \eqref{eq:moment-bochner} is the usual trace Bochner identity. When both bundles are line bundles, the two derivative terms have total
coefficient $p^2\mu^{p-1}$, as direct differentiation requires.
\end{remark}

\subsection{Differentiation of a spectral trace}

For a real $C^1$ function $q$, define its divided difference by
\[
 q^{[1]}(s,t)=
 \begin{cases}
  (q(s)-q(t))/(s-t),&s\ne t,\\
  q'(s),&s=t.
 \end{cases}
\]

We record the standard divided-difference formulas in a form suited to the $C^2$ spectral truncation used below \cite[Section~4]{Sze18}.
For completeness, we give a proof that also covers repeated eigenvalues.

\begin{lemma}
\label{lem:truncated-spectral-calculus}
Let $B$ be a $C^2$ Hermitian matrix-valued function.
At the point under consideration, suppose that
\[
 B(x)=\operatorname{diag}(\mu_1,\ldots,\mu_r).
\]
For a real-valued $C^2$ function $f$ defined on an open interval
containing the spectrum of $B$, the spectral trace
$\operatorname{tr}f(B)$ is $C^2$, and
\begin{align}
 \partial_i\partial_{\bar i}\operatorname{tr}f(B)
 ={}&\sum_{a=1}^r f'(\mu_a)
          (\partial_i\partial_{\bar i}B)_{a\bar a}
 \notag\\
 &+\sum_{a,b=1}^r (f')^{[1]}(\mu_a,\mu_b)
          |(\partial_iB)_{a\bar b}|^2.
 \label{eq:truncated-spectral-hessian}
\end{align}
For a real-valued $C^1$ function $q$ on such an interval, the map $B\mapsto q(B)$ is $C^1$ on the open set of Hermitian matrices
with spectrum in that interval, regarded as a subset of a real vector space. For a Hermitian matrix $V$, write
\[
 Dq(B)[V]=\left.\frac{d}{dt}\right|_{t=0}q(B+tV),
 \qquad t\in\mathbb R.
\]
Then
\[
 (Dq(B)[V])_{a\bar b}
 =q^{[1]}(\mu_a,\mu_b)V_{a\bar b}.
\]
The same formula holds for the complex-linear extension of the
differential, as used for $\partial_iB$.
\end{lemma}

\begin{proof}
We first consider polynomial functions. For the monomial
$q(t)=t^k$, where $k\ge1$, the product rule gives
\[
 (B+tV)^k
 =B^k+t\sum_{j=0}^{k-1}B^jVB^{k-1-j}+O(t^2).
\]
Hence
\[
 Dq(B)[V]=\sum_{j=0}^{k-1}B^jVB^{k-1-j}.
\]
At $B=\operatorname{diag}(\mu_1,\ldots,\mu_r)$,
its $(a,\bar b)$ entry is
\[
 (Dq(B)[V])_{a\bar b}
 =
 \left(\sum_{j=0}^{k-1}
 \mu_a^j\mu_b^{k-1-j}\right)V_{a\bar b}.
\]
The coefficient satisfies
\[
 \sum_{j=0}^{k-1}\mu_a^j\mu_b^{k-1-j}
 =
 \begin{cases}
 \dfrac{\mu_a^k-\mu_b^k}{\mu_a-\mu_b},
     &\mu_a\ne\mu_b,\\[6pt]
 k\mu_a^{k-1},&\mu_a=\mu_b.
 \end{cases}
\]
This is precisely $q^{[1]}(\mu_a,\mu_b)$.
Linearity proves the differential formula for all polynomials.

Taking the trace and using cyclic invariance gives
\[
\begin{aligned}
 D\operatorname{tr}(B^k)[V]
 &=\sum_{j=0}^{k-1}
   \operatorname{tr}(B^jVB^{k-1-j})
 =k\operatorname{tr}(B^{k-1}V).
\end{aligned}
\]
Consequently, for every polynomial $f$,
\[
 D\operatorname{tr}f(B)[V]
 =\operatorname{tr}(f'(B)V).
\]
By the complex-linear extension of the differential,
\[
 \partial_{\bar i}\operatorname{tr}f(B)
 =\operatorname{tr}\bigl(f'(B)\partial_{\bar i}B\bigr).
\]
Differentiating once more yields
\[
\begin{aligned}
 \partial_i\partial_{\bar i}\operatorname{tr}f(B)
 ={}&\operatorname{tr}\bigl(
 f'(B)\partial_i\partial_{\bar i}B\bigr)
 +\operatorname{tr}\bigl(
 Df'(B)[\partial_iB]\,\partial_{\bar i}B\bigr).
\end{aligned}
\]
At the diagonal matrix $B(x)$, the first term is
\[
 \sum_a f'(\mu_a)
       (\partial_i\partial_{\bar i}B)_{a\bar a}.
\]
Applying the differential formula to $f'$, the second term is
\[
 \sum_{a,b}(f')^{[1]}(\mu_a,\mu_b)
 (\partial_iB)_{a\bar b}
 (\partial_{\bar i}B)_{b\bar a}.
\]
Since $B$ is Hermitian,
\[
 (\partial_{\bar i}B)_{b\bar a}
 =\overline{(\partial_iB)_{a\bar b}}.
\]
This proves \eqref{eq:truncated-spectral-hessian}
for polynomials.

We now pass to the stated regularity classes. Choose a compact interval $J$ in the domain of the scalar function that contains the spectra of all matrices in a sufficiently small
neighborhood of $B(x)$ in the real vector space of Hermitian matrices.
For $q\in C^1$, approximate $q'$ uniformly on $J$ by polynomials and integrate, fixing the value at one point.
The resulting polynomials $q_n$ converge to $q$ in $C^1(J)$. The identity
\[
 q^{[1]}(s,t)
 =\int_0^1q'((1-\theta)t+\theta s)\,d\theta
\]
implies uniform convergence of the divided differences.

By the entrywise differential formula and unitary diagonalization, the uniform convergence of the divided differences gives uniform convergence of the differentials on this matrix neighborhood. The matrix functions also converge uniformly, since $q_n\to q$ uniformly on $J$.
The standard theorem on uniform convergence of functions and their first derivatives therefore shows that $C\mapsto q(C)$ is $C^1$, with the asserted divided-difference formula.
Its complex-linear extension satisfies the same entrywise identity.

For $f\in C^2$, applying this argument to $f$ gives
\[
 D\operatorname{tr}f(C)[V]
 =\operatorname{tr}(f'(C)V).
\]
Applying it to $f'$ shows that $C\mapsto f'(C)$ is $C^1$. Hence $C\mapsto\operatorname{tr}f(C)$ is $C^2$, and the preceding differentiation proves
\eqref{eq:truncated-spectral-hessian} in full generality. No simplicity assumption on the eigenvalues is needed.
\end{proof}

Related estimates for individual extreme eigenvalues use smooth supports
\cite[Lemma~5]{BCD17}, \cite[Lemma~3.2]{TW21}, or perturbations
\cite[Section~5]{CTW19}, \cite[Section~4]{Sze18}.
The spectral trace formula above applies directly at repeated eigenvalues.

\section{Semipositive curvature and parallel obstructions}
\label{sec:semipositive}

On a $G$-parallel holomorphic summand contained in $\ker S_b(G)$, constant enlargement of $G$ leaves the curvature contraction of the rescaled metric zero on that summand and can violate metric comparison.
We show that the eigenspaces of $G^{-1}H$ above $\kappa$ together form a $G$-parallel holomorphic summand contained in $\ker S_b(G)$.
We obtain this summand as the image of a parallel spectral truncation, without assuming regularity of the curvature kernel.

For a fixed Hermitian base metric $b$ and a Hermitian bundle metric $G$,
set
\[
 S_b(G)(\xi,\bar\eta)=b^{i\bar j}R^G_{i\bar j\xi\bar\eta},
 \qquad G(\mathscr S_G\xi,\bar\eta)=S_b(G)(\xi,\bar\eta).
\]
The fiber indices of $R^G$ are lowered with $G$. Thus $S_b(G)$ is a Hermitian form, and $\mathscr S_G$ is its $G$-self-adjoint endomorphism.
All comparisons use the same base metric $b$ and the same fiber vectors.
In particular,
\[
 S_b(tG)=tS_b(G)\qquad(t>0\text{ constant}),
\]
although the endomorphism-valued Chern curvature is unchanged by this scaling. For a K\"ahler metric $g$, one has
$S_g(g)=\operatorname{Ric}(g)$. Annihilation by $\mathscr S_G$ concerns
only the contracted curvature, not the full Chern curvature. Metrics with vanishing second Chern--Ricci curvature are studied in \cite{BP25}.

\subsection{The parallel obstruction}

\begin{proof}[Proof of Theorem~\ref{thm:excess-main}]
Let $A=G^{-1}H$, so $H(\xi,\bar\eta)=G(A\xi,\bar\eta)$.
Set $f(t)=(t-\kappa)_+^3$, where $s_+=\max\{s,0\}$, and define
\[
 F=\operatorname{tr}f(A).
\]
Let $\mathbf G,\mathbf H$ be the ordinary Gram matrices.
The Hermitian matrix
\[
 B=\mathbf G^{-1/2}\mathbf H\mathbf G^{-1/2}
\]
has the same eigenvalues as the intrinsic relative endomorphism $A$.
Lemma~\ref{lem:truncated-spectral-calculus} therefore gives
$F\in C^2(X)$.

Fix $x\in X$. Choose tangent coordinates with $b_{i\bar j}(x)=\delta_{ij}$ and a holomorphic frame normal for $G$.
A constant unitary change gives
\[
 \mathbf G(x)=I,\qquad d\mathbf G(x)=0,\qquad
 \mathbf H(x)=\operatorname{diag}(\mu_1,\ldots,\mu_r).
\]
Put $C_i=\partial_i\mathbf H(x)$ and $\rho_a=S_b(G)(e_a,\bar e_a)\ge0$.
For $S=\mathbf G^{-1/2}$, differentiating $S\mathbf G S=I$ at $x$ gives $dS=0$ and
$\partial_i\partial_{\bar i}S=-\frac12\partial_i\partial_{\bar i}\mathbf G$.
Hence, at $x$,
\[
 \partial_iB=C_i,\qquad
 \partial_i\partial_{\bar i}B
 =\partial_i\partial_{\bar i}\mathbf H
 -\frac12\bigl(
  (\partial_i\partial_{\bar i}\mathbf G)\mathbf H
  +\mathbf H(\partial_i\partial_{\bar i}\mathbf G)\bigr).
\]
The Chern curvature formula, with its fiber index lowered, reads
\[
 R^H_{i\bar i}=-\partial_i\partial_{\bar i}\mathbf H
                   +C_i\mathbf H^{-1}C_i^*,
 \qquad
 R^G_{i\bar i}=-\partial_i\partial_{\bar i}\mathbf G.
\]
Consequently, with
$Q=\sum_i C_i\mathbf H^{-1}C_i^*\ge0$, the spectral Hessian formula becomes
\begin{align}
 \Delta_bF={}&
 \sum_a f'(\mu_a)
 \bigl\{Q_{a\bar a}+\mu_a\rho_a
                      -S_b(H)(e_a,\bar e_a)\bigr\}
 \notag\\
 &+\sum_{i,a,b}(f')^{[1]}(\mu_a,\mu_b)
                    |(C_i)_{a\bar b}|^2.
 \label{eq:truncated-comparison-exact}
\end{align}
Here $\Delta_b=\sum b^{i\bar j}\partial_i\partial_{\bar j}$ on functions. Both $f'$ and $(f')^{[1]}$ are nonnegative.
Therefore the curvature comparison gives
\begin{align}
 \Delta_bF
 \ge{}&\sum_a f'(\mu_a)Q_{a\bar a}
  +\sum_{i,a,b}(f')^{[1]}(\mu_a,\mu_b)
                    |(C_i)_{a\bar b}|^2
 \notag\\
 &+3\sum_a\rho_a(\mu_a-\kappa)_+^3
 \ge0.
 \label{eq:truncated-comparison-positive}
\end{align}
Compactness and the strong maximum principle imply that $F$ is constant. Thus every nonnegative term in
\eqref{eq:truncated-comparison-positive} vanishes at every point.

The divided-difference term controls every derivative component involving the excess spectrum. By symmetry in $s,t$ and the formulas
\[
 (f')^{[1]}(s,t)=
 \begin{cases}
  3(s+t-2\kappa),&s,t>\kappa,\\[2pt]
  \dfrac{3(s-\kappa)^2}{s-t},&s>\kappa\ge t,
 \end{cases}
\]
this coefficient is positive whenever $\max\{s,t\}>\kappa$.
Consequently,
\[
 (C_i)_{a\bar b}=0
 \qquad\text{if }\max\{\mu_a,\mu_b\}>\kappa.
\]

Set $T=f(A)=(A-\kappa I)_+^3$. Orthonormalizing the chosen normal frame gives a smooth $G$-unitary frame with the same value and vanishing connection form at $x$.
The matrix of $A$ is Hermitian in this frame, so Lemma~\ref{lem:truncated-spectral-calculus} gives $T\in C^1$ and
\[
 G((\nabla_i^GT)e_a,\bar e_b)
 =f^{[1]}(\mu_a,\mu_b)G((\nabla_i^GA)e_a,\bar e_b)
 =f^{[1]}(\mu_a,\mu_b)(C_i)_{a\bar b}.
\]
If both eigenvalues are at most $\kappa$, then $f^{[1]}(\mu_a,\mu_b)=0$; otherwise $(C_i)_{a\bar b}=0$.
Thus $\nabla_i^GT=0$. Metric compatibility and the $G$-self-adjointness of $T$ also give
$\nabla_{\bar i}^GT=0$, so $\nabla^GT=0$.

Parallel transport conjugates $T$ between fibers. It follows that $T$ is smooth and has constant eigenvalues and rank on connected $X$.
Its image is exactly $E_+$, which is therefore a smooth parallel subbundle. Metric compatibility gives the parallel splitting
\[
 E=E_+\oplus E_+^{\perp_G}.
\]
This splitting is holomorphic because $(\nabla^G)^{0,1}=\bar\partial_E$.
The identity $A|_{E_+}=\kappa I+(T|_{E_+})^{1/3}$ also shows that $A|_{E_+}$ is parallel.

Finally, the last sum in \eqref{eq:truncated-comparison-positive} gives $\rho_a=0$
whenever $\mu_a>\kappa$. Semipositivity then implies
\[
 |S_b(G)(e_a,\bar v)|^2
 \le S_b(G)(e_a,\bar e_a)S_b(G)(v,\bar v)=0
\]
for every $v\in E_x$ and each such $e_a$. These vectors span $(E_+)_x$, so $\mathscr S_G|_{E_+}=0$.
If there is no nonzero parallel subbundle annihilated by $\mathscr S_G$, then $E_+=0$, equivalently $H\le\kappa G$.
\end{proof}

\begin{remark}\label{rem:truncation-power}
The cubic is the lowest integer power for which the truncation is
$C^2$ at $\kappa$. Any real exponent $\alpha>2$ works in the same proof:
for $f_\alpha(t)=(t-\kappa)_+^\alpha$,
\[
 f_\alpha'(t)(t-\kappa)=\alpha(t-\kappa)_+^\alpha,
 \qquad f_\alpha''(t)>0\quad(t>\kappa).
\]
The divided difference of $f_\alpha'$ is positive whenever at least
one argument is greater than $\kappa$, so the equality terms again give parallelism. The quadratic truncation is only $C^{1,1}$ and
does not directly fit the classical $C^2$ calculation used here.
\end{remark}

\subsection{A criterion for comparison}

The obstruction in Theorem~\ref{thm:excess-main} is also necessary if comparison is required for every metric $H$.

\begin{corollary}
\label{cor:universal-comparison}
Fix a compact connected Hermitian manifold $(X,b)$, a Hermitian holomorphic bundle $(E,G)$ with $S_b(G)\ge0$, and a constant
$\kappa>0$. The following are equivalent:
\begin{enumerate}
\item[(i)] Every Hermitian metric $H$ on $E$ satisfying $S_b(H)\le\kappa S_b(G)$ also satisfies $H\le\kappa G$.
\item[(ii)] There is no nonzero $G$-parallel subbundle of $E$ annihilated by $\mathscr S_G$.
\end{enumerate}
\end{corollary}

\begin{proof}
The implication (ii)$\Rightarrow$(i) is Theorem~\ref{thm:excess-main}. Conversely, suppose such a subbundle $V$ exists. Its $G$-orthogonal projection $P$ is parallel.
For $t>\kappa$, set
\[
 A_t=tP+\kappa(I-P),\qquad
 H_t(\xi,\bar\eta)=G(A_t\xi,\bar\eta).
\]
Since $A_t$ is parallel, the Chern connection of $G$ is compatible with $H_t$ and has $(0,1)$ part $\bar\partial_E$. By uniqueness, it is also the Chern connection of $H_t$. Since
$\mathscr S_G$ is self-adjoint and annihilates $V$, one has $A_t\mathscr S_G=\kappa\mathscr S_G$. Lowering the fiber index with
$H_t$ therefore gives
\[
 S_b(H_t)(\xi,\bar\eta)
 =G(A_t\mathscr S_G\xi,\bar\eta)
 =\kappa S_b(G)(\xi,\bar\eta).
\]
But $H_t=tG$ on $V$, so (i) fails.
\end{proof}

A nontrivial curvature kernel need not contain a parallel subbundle. The following quotient metric illustrates this distinction on a curve, with a one-dimensional kernel at every point.

\begin{example}
\label{ex:veronese-comparison}
Let $\nu:\mathbb{CP}^1\to\mathbb{CP}^2$ be the Veronese embedding
\[
 \nu([s:t])=[s^2:st:t^2],
\]
and let $(E,G)=\nu^*Q$, where $Q$ is the universal quotient bundle with the metric induced by the standard metric on $\mathbb C^3$. For every Hermitian base metric $b$ on
$\mathbb{CP}^1$, the form $S_b(G)$ is semipositive of complex rank one everywhere. Nevertheless, every Hermitian metric $H$ on $E$ satisfying $S_b(H)\le\kappa S_b(G)$, where $\kappa>0$, obeys $H\le\kappa G$.

To verify these assertions, pull back the tautological sequence:
\[
 0\longrightarrow\mathcal O(-2)
 \xrightarrow{(s^2,st,t^2)^T}\mathcal O^{\oplus3}
 \longrightarrow E\longrightarrow0.
\]
Here $\mathcal O(j)$ denotes the line bundle of degree $j$ on $\mathbb{CP}^1$. In the dual sequence, the kernel of the row $(s^2,st,t^2)$ is the image of
\[
 \mathcal O(-1)^{\oplus2}\longrightarrow\mathcal O^{\oplus3},
 \qquad
 \begin{pmatrix}
 -t&0\\ s&-t\\ 0&s
 \end{pmatrix}.
\]
The columns lie in the kernel and are linearly independent whenever $(s,t)\ne(0,0)$. Hence $E^*\simeq\mathcal O(-1)^{\oplus2}$.

Let $L=\mathcal O(-2)$ be the tautological subbundle in the displayed sequence, and let $\beta_v:L_x\to E_x$ be its second fundamental
form in a tangent direction $v$. The curvature formula for a quotient of a flat Hermitian bundle \cite{Gr69} gives
\[
 \mathcal R^G(v,\bar v)=\beta_v\beta_v^*\ge0,
\]
where $G(\mathcal R^G(v,\bar v)\xi,\bar\eta)=R^G(v,\bar v,\xi,\bar\eta)$. The map $\beta_v$ has rank one
for $v\ne0$: in the chart $s\ne0$, the derivative of $(1,z,z^2)$ is not a multiple of $(1,z,z^2)$; the same argument
applies in the other chart. Since the base has complex dimension one, contracting this curvature gives $\operatorname{rank}_{\mathbb C}\mathscr S_G=1$ at every point.

Suppose a nonzero $G$-parallel subbundle $V$ were annihilated by $\mathscr S_G$. Since $V$ is parallel, its induced curvature is the restriction of the curvature of $G$. On a complex curve this curvature is determined by its contraction, which vanishes. Hence $V$ is flat.
Since $\mathbb{CP}^1$ is simply connected, it gives a holomorphic trivialization $V\simeq\mathcal O^{\oplus a}$, where $a=\operatorname{rank}_{\mathbb C}V>0$.
The parallel orthogonal splitting makes $V$ a holomorphic direct summand of $E$. Composing the projection $E\to V$ with a nonzero
coordinate projection $V\to\mathcal O$ yields a nonzero section of $E^*$. This contradicts
$H^0(\mathbb{CP}^1,E^*)=H^0(\mathbb{CP}^1,\mathcal O(-1))^{\oplus2}=0$.
Corollary~\ref{cor:universal-comparison} now gives the asserted
comparison.
\end{example}

For this fixed metric $G$, the positive-definite contracted curvature hypothesis of \cite{WYY26} fails everywhere. The strict directional positivity required in \cite[Theorem~1.8]{XYY} also fails: on a curve, a nonzero vector in $\ker\mathscr S_G$ has zero curvature in every
base direction. Thus Theorem~\ref{thm:excess-main} applies beyond these positivity assumptions. The comparison concerns this prescribed quotient metric $G$.

\subsection{Positive comparison}

Strict positivity excludes the obstruction. Positive definiteness at one point already suffices.

\Needspace{10\baselineskip}
\begin{corollary}\label{thm:metric-comparison}
Let $(X,b)$ be compact connected Hermitian, and let $G,H$ be Hermitian metrics on a holomorphic vector bundle $E\to X$. Suppose
\[
 S_b(G)\ge0,\qquad S_b(H)\le\kappa S_b(G),\qquad \kappa>0.
\]
If $S_b(G)$ is positive definite at one point, then $H\le\kappa G$.
\end{corollary}

\begin{proof}
The parallel subbundle $E_+$ in Theorem~\ref{thm:excess-main} has constant rank and is annihilated by $\mathscr S_G$. It is zero at the point of positive
definiteness, hence zero everywhere.
\end{proof}

\begin{remark}\label{rem:metric-xyy}
The case $S_b(G)>0$ follows from \cite[Theorem~1.8]{XYY}. Fix $x\in X$, put $d=\dim_{\mathbb C}X$,
and choose $0\ne\xi\in E_x$. Set $P(v,\bar v)=R^G(v,\bar v,\xi,\bar\xi)$ and $Q(v,\bar v)=R^H(v,\bar v,\xi,\bar\xi)$.
Their traces satisfy $\operatorname{tr}P>0$ and $\operatorname{tr}Q\le\kappa\operatorname{tr}P$.
Applying the discrete Fourier transform to a $b$-unitary eigenbasis of $P$ gives a $b$-unitary basis $w_1,\ldots,w_d$ with
$P(w_j,\bar w_j)=d^{-1}\operatorname{tr}_bP>0$ for every $j$. Since $\sum_j(Q-\kappa P)(w_j,\bar w_j)\le0$, at least one $w_j$
also satisfies $Q(w_j,\bar w_j)\le\kappa P(w_j,\bar w_j)$. This is the common direction required by Xiong--Yang--Yau.
On a compact K\"ahler base, the contracted comparison is stated in \cite[Theorem~1.2 and Corollary~3.4]{WYY26}.
Positive definiteness at one point is covered by Corollary~\ref{thm:metric-comparison}.
On a compact complex manifold, Lee--Li \cite[Theorem~1.1]{LL20} prove that a tangent-bundle metric whose Chern curvature contraction
is nonnegative everywhere and positive definite at one point can be conformally deformed to make that contraction positive definite
everywhere, with the base metric used for contraction fixed.
An alternative proof for $S_b(G)>0$ using weighted moments is given in Remark~\ref{rem:weighted-positive}.
\end{remark}

\subsection{Tangent metrics and a product obstruction}

\begin{corollary}
\label{cor:semipositive-tangent-comparison}
Let $(M,g)$ be a compact connected K\"ahler manifold with
$\operatorname{Ric}_g\ge0$, and let $h$ be a Hermitian metric on $M$.
Assume, for a constant $\kappa>0$, that
\[
 \operatorname{tr}_{\omega_g}R^h
 \le\kappa\operatorname{Ric}_g.
\]
If $T^{1,0}M$ has no nonzero $g$-parallel complex subbundle annihilated
by the Ricci endomorphism, then $h\le\kappa g$. In particular, this
holds if $g$ is locally irreducible and
$\operatorname{Ric}_g\not\equiv0$.
\end{corollary}

\begin{proof}
Apply Theorem~\ref{thm:excess-main} with
$E=T^{1,0}M$, $b=G=g$, and $H=h$. The K\"ahler curvature symmetries
give $S_g(g)=\operatorname{Ric}_g$. The Chern connection is the restriction
of the complexified Levi--Civita connection. If $g$ is locally
irreducible, a nonzero parallel complex subbundle must be the whole
tangent bundle: its underlying real distribution and its orthogonal
complement would otherwise give a local Riemannian product. The
condition $\operatorname{Ric}_g\not\equiv0$ rules out the whole tangent
bundle as a Ricci-annihilated subbundle.
\end{proof}

The positive case also covers a complete K\"ahler metric $g$ satisfying $\operatorname{Ric}_g\ge cg$ for some $c>0$.
The Bonnet--Myers theorem \cite{Mye41} makes $M$ compact, and strict
Ricci positivity excludes the parallel obstruction.
Thus every Hermitian metric $h$ with $S_g(h)\le\kappa\operatorname{Ric}_g$, where $\kappa>0$, satisfies
$h\le\kappa g$ by Corollary~\ref{cor:semipositive-tangent-comparison}.

\Needspace{10\baselineskip}
\begin{corollary}
\label{cor:metric-rigidity}\label{cor:common-volume}
Under the assumptions of
Corollary~\ref{cor:semipositive-tangent-comparison}, let
$m=\dim_{\mathbb C}M$. Then
\[
 \operatorname{Vol}(M,h)\le\kappa^m\operatorname{Vol}(M,g),
\]
with equality if and only if $h=\kappa g$.
\end{corollary}

\begin{proof}
The comparison $h\le\kappa g$ puts all eigenvalues of $g^{-1}h$ in
$(0,\kappa]$. Their product is the volume density ratio, so the volume bound follows. Equality of the integrals forces this continuous
product to equal $\kappa^m$ everywhere. Each eigenvalue must then equal $\kappa$, which gives $h=\kappa g$; the converse is immediate.
\end{proof}

The contracted curvature hypothesis in these tangent-metric comparisons
follows from $R^h(v,\bar v,w,\bar w)\le\kappa R^g(v,\bar v,w,\bar w)$
for all $v,w\in T^{1,0}M$, by tracing over $v$.
The volume equality criterion uses only the metric comparison, so in
the strictly positive case it also follows from \cite{WYY26,XYY}.

\begin{example}
\label{ex:ricci-flat-parallel-obstruction}
Let $M=\mathbb{CP}^1\times T$, where $T=\mathbb C/\Lambda$ is an
elliptic curve, $g_1$ is a Fubini--Study metric, and $g_0$ is a flat
metric on $T$. Set
\[
 g=g_1\oplus g_0,
 \qquad h=\kappa g_1\oplus c g_0,
 \qquad c>\kappa>0.
\]
Constant rescaling leaves each Chern connection unchanged and scales
the corresponding curvature tensor with its fiber index lowered.
The flat factor contributes no curvature. Consequently,
\[
 \operatorname{tr}_{\omega_g}R^h
 =\kappa\operatorname{Ric}_g,
 \qquad \operatorname{Ric}_g\ge0,
\]
but $h\nleq\kappa g$. Here the subbundle $E_+$ in
Theorem~\ref{thm:excess-main} is exactly the tangent line
of the elliptic factor. Thus semipositive Ricci curvature alone does not imply the metric
comparison. The excess spectrum is confined to the Ricci-flat
parallel factor.
\end{example}

The exceptional factor is Ricci-flat, not necessarily flat. The same
construction works with any compact Ricci-flat K\"ahler factor (e.g., a $K3$ surface) in
place of $T$: its contracted curvature is zero, even when its full
curvature is nonzero.

\section{Schwarz estimates on complete bases}\label{sec:complete}

The main analytic step is to bound the Schatten moments $U_p$
without an a priori upper bound. The lower bound on $\operatorname{tr}_{\omega_g}R^{E_1}$
and H\"older's inequality give a superlinear differential inequality,
to which the following maximum principle applies.

\subsection{A nonlinear maximum principle}

We use a standard consequence of the Omori--Yau maximum principle \cite{Omori67,Yau75}. Omori's Hessian principle assumes a sectional
curvature lower bound; the Ricci lower bound here suffices for Yau's
Laplacian version. We use the complex Laplacian and gradient convention
of Section~\ref{sec:preliminaries}.

\Needspace{9\baselineskip}
\begin{lemma}\label{lem:nonlinear-max}
Let $(M,g)$ be a complete K\"ahler manifold with Ricci curvature bounded
from below.  If a nonnegative $C^2$ function $u$ satisfies
\[
 \Delta_g u\ge a u^{1+\delta}-bu
\]
for constants $a,\delta>0$ and $b\ge0$, then
\[
 \sup_M u\le (b/a)^{1/\delta}.
\]
In particular, if $b=0$, then $u\equiv0$.
\end{lemma}

\begin{proof}
Set $\alpha=\delta/2$ and $F=(u+1)^{-\alpha}$.  The function $F$ is
positive and bounded from below.  By Yau's maximum principle \cite{Yau75}, for every $\varepsilon>0$ there is $x_\varepsilon\in M$
such that
\[
 F(x_\varepsilon)<\inf_MF+\varepsilon,\qquad
 |\partial F|_g(x_\varepsilon)<\varepsilon,\qquad
 \Delta_gF(x_\varepsilon)>-\varepsilon.
\]
Direct differentiation gives
\[
 \Delta_gF=-\alpha(u+1)^{-\alpha-1}\Delta_gu
 +\alpha(\alpha+1)(u+1)^{-\alpha-2}|\partial u|_g^2
\]
and
\[
 |\partial u|_g^2=\alpha^{-2}(u+1)^{2\alpha+2}|\partial F|_g^2.
\]
Since $\delta=2\alpha$, the two Omori--Yau inequalities at
$x_\varepsilon$ therefore imply
\begin{equation}\label{eq:nonlinear-max}
 \alpha(u+1)^{-1-\delta}\Delta_gu
 \le \frac{\alpha+1}{\alpha}\varepsilon^2
      +\varepsilon(u+1)^{-\alpha}
 \quad\text{at }x_\varepsilon.
\end{equation}
If $\sup_Mu=\infty$, then $u(x_\varepsilon)\to\infty$ along a sequence
$\varepsilon\to0$.  Substitution of
$\Delta_gu\ge u(au^\delta-b)$ into \eqref{eq:nonlinear-max} gives a lower bound for its
left-hand side tending to $\alpha a$, while the right-hand side tends
to zero, a contradiction.  Thus $S=\sup_Mu<\infty$, and the same choice of points satisfies
$u(x_\varepsilon)\to S$. If $S>0$, letting $\varepsilon\to0$ in \eqref{eq:nonlinear-max} gives $aS^\delta-b\le0$.
The case $S=0$ is immediate.
\end{proof}

\subsection{The complete estimate}

\begin{proof}[Proof of Theorem~\ref{thm:complete-main}]
If $\varphi\equiv0$, the conclusion is immediate. We may therefore
assume $q\ge1$. Fix an integer $p\ge2$ and set
\[
 U_p=\operatorname{tr}(A^p)
     =\sum_{a=1}^{q}\mu_a^p.
\]
Tracing the inequality in Proposition~\ref{prop:moment-bochner} with
respect to $g$ and using the two curvature assumptions gives
\begin{equation}\label{eq:complete-up}
 \Delta_gU_p\ge pK_1\sum_{a=1}^{q}\mu_a^{p+1}
                 -pK_2E_fU_p .
\end{equation}
H\"older's inequality gives
\[
 U_p=\sum_{a=1}^{q}(\mu_a^{p+1})^{p/(p+1)}\cdot1
 \le
 \left(\sum_{a=1}^{q}\mu_a^{p+1}\right)^{p/(p+1)}
 q^{1/(p+1)}.
\]
Raising both sides to the power $(p+1)/p$ and rearranging, we obtain
\begin{equation}\label{eq:holder-moment}
 \sum_{a=1}^{q}\mu_a^{p+1}
 \ge q^{-1/p}U_p^{1+1/p}.
\end{equation}
Consequently
\[
 \Delta_gU_p\ge
 pK_1q^{-1/p}U_p^{1+1/p}-pK_2E_fU_p.
\]
Lemma~\ref{lem:nonlinear-max}, with $\delta=1/p$, yields
\[
 \sup_M U_p\le q\left(\frac{K_2E_f}{K_1}\right)^p.
\]
Since $\mu_1\le U_p^{1/p}$, we obtain
\[
 \mu_1\le q^{1/p}\frac{K_2E_f}{K_1}.
\]
Letting $p\to\infty$ proves the asserted tensor inequality.
\end{proof}

For general bundle maps, the assumption $E_f<\infty$ controls
the energy density of the base map $f$, without imposing a prior
bound on $\varphi$. When $\varphi=df$, however, this assumption
already gives a bound on the bundle map itself.

\begin{remark}\label{rem:sharp-complete}
The constant in Theorem~\ref{thm:complete-main} cannot in general be
improved.  Let $M=N=\mathbb{CP}^1$ with the same Fubini--Study metric,
let $f$ be the identity, and equip $L=\mathcal O(1)$ with
its standard Hermitian metric $G$. Its curvature satisfies
\[
 R^G(w,\bar w,\xi,\bar\xi)
 =\kappa|w|_g^2|\xi|_G^2
\]
for a positive constant $\kappa$.
Take $E_1=E_2=L$, let $H=cG$ for a constant $c>0$, and let
$\varphi$ be the identity bundle map.  Then $E_f=1$,
$\varphi^*H=cG$, and constant rescaling does not change the Chern
connection.  The curvature assumptions hold with
$K_1=\kappa/c$ and $K_2=\kappa$.  Hence the conclusion reads
$\varphi^*H\le cG$, and equality holds everywhere.
Taking $r$ copies of $L$ gives $\operatorname{tr}_G\varphi^*H=rc$
with the same curvature constants. Thus a rank factor is unavoidable
in the corresponding trace bound.

Equality also occurs on a complete noncompact base. Take
$M=\mathbb{CP}^1\times\mathbb C$ with the product of the same
Fubini--Study metric and the Euclidean metric, let $N=\mathbb{CP}^1$,
and let $f$ be the projection onto the first factor. Set
$E_1=f^*L$ with metric $f^*G$, $E_2=L$ with metric $H=cG$, and let
$\varphi$ be the canonical bundle map covering $f$. The product metric
has nonnegative Ricci curvature, $E_f=1$, and
$\operatorname{tr}_{\omega_g}R^{E_1}=\kappa f^*G$.
Thus the same constants $K_1=\kappa/c$, $K_2=\kappa$ give equality
$\varphi^*H=c f^*G$ in Theorem~\ref{thm:complete-main}.
\end{remark}

\subsection{Compact variants of the moment argument}

We next compare $\operatorname{tr}_{\omega_g}R^{E_1}$ with $G$ and bound the target curvature on $\varphi(\xi)$ in terms of $|\xi|_G^2$.
The resulting moment inequality has a positive term linear in $U_p$, so we argue at a maximum point on a compact base.

\Needspace{14\baselineskip}
\begin{prop}\label{prop:compact-companion}
Let $(M,g)$ be compact Hermitian and let $(N,h)$ be Hermitian.  Suppose
that, for constants $K_1,K_2>0$,
\[
 \operatorname{tr}_{\omega_g}R^{E_1}\ge K_1G
\]
and
\[
 R^{E_2}(df_x(v),\overline{df_x(v)},
          \varphi_x(\xi),\overline{\varphi_x(\xi)})
 \le K_2|v|_g^2|\xi|_G^2
\]
for every $x\in M$, $v\in T_x^{1,0}M$, and $\xi\in(E_1)_x$.  Then
\[
 \varphi^*H\le\frac{r_fK_2}{K_1}G.
\]
\end{prop}

\begin{proof}
Let $p\ge2$ and choose a maximum point $z_p$ of $U_p$. Use a $g$-unitary basis at $z_p$ whose last $m-\operatorname{rank}_{\mathbb C}df_{z_p}$ vectors lie in $\ker df_{z_p}$.
If $\mu_a>0$, the relation $\varphi(e_a)=\lambda_as_a$ bounds the target curvature by $K_2/\mu_a$
in each remaining basis direction. The terms with $\mu_a=0$ vanish. Proposition~\ref{prop:moment-bochner} therefore gives
\begin{equation}\label{eq:compact-companion-up}
 0\ge\Delta_gU_p(z_p)
 \ge pK_1U_p(z_p)-pr_fK_2U_{p-1}(z_p).
\end{equation}
H\"older's inequality yields
$U_{p-1}\le q^{1/p}U_p^{(p-1)/p}$, and hence
\[
 \max_M U_p^{1/p}\le q^{1/p}\frac{r_fK_2}{K_1}.
\]
The case $U_p(z_p)=0$ is immediate. Since $\mu_1\le U_p^{1/p}$, the limit $p\to\infty$ proves the result.
\end{proof}

\begin{remark}\label{rem:companion-noncompact}
Proposition~\ref{prop:compact-companion} is false on a general complete noncompact base.  Take $M=N=\mathbb C$ with the Euclidean metric and
$f=\operatorname{id}$.  Let $E_1$ and $E_2$ be trivial line bundles with global frames $e_1,e_2$.  Set
\[
 G(e_1,\bar e_1)=e^{-|z|^2},\qquad H(e_2,\bar e_2)=1,
 \qquad \varphi(e_1)=e_2.
\]
With our sign convention, $\operatorname{tr}_{\omega_g}R^{E_1}=G$ and $E_2$ is flat.  Thus the hypotheses hold with $K_1=K_2=1$ and $r_f=1$.
However,
\[
 \frac{\varphi^*H}{G}=e^{|z|^2}
\]
is unbounded.  In the proof above, the positive term is only linear in $U_p$.  The superlinear term in \eqref{eq:holder-moment}, which is the
key to Theorem~\ref{thm:complete-main}, is absent.
\end{remark}

\begin{remark}
\label{rem:weighted-positive}
When $S_b(G)>0$, the moment method also recovers the comparison in Remark~\ref{rem:metric-xyy}.
The sharp constant is retained by using the same curvature weights at consecutive exponents.
Let $G,H$ be metrics on $E$ over a compact Hermitian base $(X,b)$ with $S_b(G)>0$ and $S_b(H)\le\kappa S_b(G)$, where $\kappa>0$.
Choose uniform constants $c,C>0$ with $cG\le S_b(G)\le CG$. For an integer $p\ge2$, take a maximum point $z_p$ of $U_p=\operatorname{tr}((G^{-1}H)^p)$. Choose there a $G$-unitary
eigenbasis $e_a$ of $G^{-1}H$, with eigenvalues $\mu_a>0$, and set
\[
 \rho_a=S_b(G)(e_a,\bar e_a),\qquad
 W_j=\sum_{a=1}^{r}\rho_a \mu_a^j,\qquad r=\operatorname{rank}_{\mathbb C}E.
\]
The target frame $s_a=\mu_a^{-1/2}e_a$ is $H$-unitary. Applying the moment formula to the identity bundle map gives
\begin{equation}\label{eq:common-weight-moment}
 0\ge\Delta_bU_p(z_p)
 \ge p\sum_a\rho_a \mu_a^{p-1}(\mu_a-\kappa)
 =p(W_p-\kappa W_{p-1}).
\end{equation}
Weighted H\"older gives
$W_{p-1}\le W_p^{(p-1)/p}W_0^{1/p}$. Thus
\[
 c\max_XU_p\le W_p\le\kappa^pW_0\le rC\kappa^p.
\]
Taking $p$-th roots and letting $p\to\infty$ gives $H\le\kappa G$. The common weights place $C/c$ inside the $p$-th root, preserving the comparison constant $\kappa$.
\end{remark}

\section{Partial spectral estimates on compact bases}\label{sec:compact}

Throughout this section $M$ is compact and $1\le k\le q$. We write $k'=\min\{k,r_\varphi\}$ and retain $E_f=\sup_M\operatorname{tr}_{\omega_g}f^*\omega_h$.
We apply Ni's local-support method for partial sums and exterior-power norms \cite[Propositions~2.1--2.2 and~3.1]{Ni19} to general bundle maps.

\subsection{Local supports and Bochner formulas}

Choose the normal frames \eqref{BJ} at a fixed point $x_0$. On a neighborhood $B_{x_0}$, let
\[
 G_{(k)}=\bigl(G(e_a,\bar e_b)\bigr)_{a,b\le k},
 \qquad
 B_{(k)}=\bigl((\varphi^*H)(e_a,\bar e_b)\bigr)_{a,b\le k}
\]
be the Gram matrices on the holomorphic span of
$e_1,\ldots,e_k$. Define
\begin{equation}\label{U1}
 \mathcal U_k=\operatorname{tr}(G_{(k)}^{-1}B_{(k)}).
\end{equation}
Then
\begin{equation}\label{UKP}
 \mathcal U_k(x_0)=\sigma_k(x_0),\qquad
 \mathcal U_k\le\sigma_k\quad\text{on }B_{x_0}.
\end{equation}
Similarly, set
\begin{equation}\label{W1}
 \mathcal W_k=\frac{\det B_{(k)}}{\det G_{(k)}}.
\end{equation}
We have
\begin{equation}\label{WKP}
 \mathcal W_k(x_0)=\pi_k(x_0),\qquad
 \mathcal W_k\le\pi_k\quad\text{on }B_{x_0}.
\end{equation}

The inequalities follow from Ni's fiberwise linear-algebra argument \cite[Propositions~2.1--2.2]{Ni19}, which applies equally to
holomorphic bundle maps. Equality holds at $x_0$ because $e_1(x_0),\ldots,e_k(x_0)$ are the first $k$ singular vectors.
Both supports are smooth on $B_{x_0}$, since the Gram matrices are smooth and $G_{(k)}$ is positive definite.
Write $\varphi^\gamma_{ai}=(\nabla_i\varphi)_a^\gamma$ in the frames \eqref{BJ}. Repeated coordinate indices are summed.

\begin{prop}\label{prop:local-bochner}
Let $\mathcal U_k$ and $\mathcal W_k$ be defined by \eqref{U1} and \eqref{W1}.  At the supporting point $x_0$ the following formulas hold.
\begin{enumerate}
\item In holomorphic coordinates $(z^1,\ldots,z^m)$,
\begin{equation}
\begin{aligned}
\partial\bar{\partial}\mathcal U_k&=\sum_{c=1}^{r_2}\sum_{a=1}^{k}\varphi^c_{ai}\overline{\varphi^c_{aj}}dz^i\wedge d\bar{z}^j \\
&\ +\sum_{a=1}^{k}\mu_a (R^{E_1}_{i\bar{j}a\bar{a}}-f^{\alpha}_{i}\overline{f^{\beta}_{j}}R^{E_2}_{\alpha{\bar{\beta}}a\bar{a}})dz^i\wedge d\bar{z}^j.
\end{aligned}
\end{equation}
Moreover, for any $v\in T^{1,0}_{x_0}M$,
\begin{equation} \label{Ukk}
\partial\bar{\partial}\mathcal U_k(v,\bar{v})\ge \sum_{a=1}^{k}\mu_a (R^{E_1}(v,\bar{v},e_a,\bar{e}_a)-R^{E_2}(f_{\ast}v,\overline{f_{\ast}v},s_a,\bar{s}_a)).
\end{equation}
\item If $\mathcal W_k(x_0)\neq0$, then under local coordinates,
\begin{equation}
\begin{aligned}
\partial\bar{\partial}\log \mathcal W_k&=\sum_{a=1}^{k}\sum_{\gamma=k+1}^{r_2}\frac{1}{\mu_a}\varphi^{\gamma}_{ai}\overline{\varphi^{\gamma}_{aj}}dz^i\wedge d\bar{z}^j \\
&\ +\sum_{a=1}^{k}R^{E_1}_{i\bar{j}a\bar{a}}dz^i\wedge d\bar{z}^j-\sum_{a=1}^{k}f^{\alpha}_{i}\overline{f^{\beta}_{j}}R^{E_2}_{\alpha{\bar{\beta}}a\bar{a}}dz^i\wedge d\bar{z}^j
\end{aligned}
\end{equation}
and for $v\in T^{1,0}_{x_0}M$,
\begin{equation}
\begin{aligned}
\partial\bar{\partial}\log \mathcal W_k(v,\bar{v})&\ge\sum_{a=1}^{k}R^{E_1}(v,\bar{v},e_a,\bar{e}_a)-\sum_{a=1}^{k}R^{E_2}(f_*v,\overline{f_*v},s_a,\bar{s}_a).
\end{aligned}
\end{equation}
\end{enumerate}
\end{prop}

\begin{proof}
Put $B=\varphi^*H$ and write $\varphi^\gamma_{ai}=(\nabla_i\varphi)_a^{\gamma}$ at $x_0$.
All calculations below are made at $x_0$ in the frames \eqref{BJ}; repeated coordinate indices are summed.  In particular,
\[
 G_{a\bar b}=\delta_{ab},\qquad
 B_{a\bar b}=\mu_a\delta_{ab},\qquad
 \partial G_{a\bar b}=0
 \qquad(1\le a,b\le k).
\]
Since the first derivatives of $G_{(k)}$ vanish at $x_0$, differentiating its inverse gives
\[
 \partial_i\partial_{\bar j}(G_{(k)}^{-1})
 =-\partial_i\partial_{\bar j}G_{(k)}
 =(R^{E_1}_{i\bar ja\bar b})_{a,b\le k}.
\]
The product rule in the definition \eqref{U1} therefore yields
\begin{equation}\label{eq:Uk-proof-split}
 \partial_i\partial_{\bar j}\mathcal U_k
 =\sum_{a=1}^{k}\mu_aR^{E_1}_{i\bar ja\bar a}
  +\sum_{a=1}^{k}\partial_i\partial_{\bar j}B_{a\bar a}.
\end{equation}

Substituting \eqref{eq:B-second-derivative} into \eqref{eq:Uk-proof-split} gives the formula for $\mathcal U_k$.
Its derivative term evaluates on $(v,\bar v)$ to $\sum_{a=1}^{k}\sum_{\gamma=1}^{r_2}|\varphi^\gamma_{ai}v^i|^2\ge0$, yielding \eqref{Ukk}.

We next assume $\mathcal W_k(x_0)>0$, so $\mu_a>0$ for $1\le a\le k$. Let $C=(B_{a\bar b})_{1\le a,b\le k}$, which is positive definite after shrinking the neighborhood of $x_0$.
For a positive-definite Hermitian matrix-valued function $C$,
\begin{equation}\label{eq:logdet-matrix}
 \partial_i\partial_{\bar j}\log\det C
 =\operatorname{tr}(C^{-1}\partial_i\partial_{\bar j}C)
 -\operatorname{tr}(C^{-1}\partial_iC\,
                    C^{-1}\partial_{\bar j}C).
\end{equation}
Using \eqref{eq:B-second-derivative} and \eqref{eq:B-first-derivatives} in \eqref{eq:logdet-matrix}, with
$C^{-1}=\operatorname{diag}(\mu_1^{-1},\ldots,\mu_k^{-1})$, gives
\begin{align*}
 \partial_i\partial_{\bar j}\log\det C
={}&\sum_{a=1}^{k}\sum_{\gamma=1}^{r_2}
     \frac{\varphi^\gamma_{ai}
     \overline{\varphi^\gamma_{aj}}}{\mu_a}
 -\sum_{a,b=1}^{k}
     \frac{\varphi^b_{ai}\overline{\varphi^b_{aj}}}{\mu_a}\\
 &-\sum_{a=1}^{k}f_i^\alpha\overline{f_j^\beta}
     R^{E_2}_{\alpha\bar\beta a\bar a}\\
={}&\sum_{a=1}^{k}\sum_{\gamma=k+1}^{r_2}
     \frac{\varphi^\gamma_{ai}
     \overline{\varphi^\gamma_{aj}}}{\mu_a}
 -\sum_{a=1}^{k}f_i^\alpha\overline{f_j^\beta}
     R^{E_2}_{\alpha\bar\beta a\bar a}.
\end{align*}
Applying \eqref{eq:logdet-matrix} to the $k\times k$ block of $G$ is even simpler, since its first derivatives vanish:
\[
 \partial_i\partial_{\bar j}\log\det(G_{a\bar b})_{a,b\le k}
 =-\sum_{a=1}^{k}R^{E_1}_{i\bar ja\bar a}.
\]
Taking the difference proves the displayed formula for $\partial\bar\partial\log \mathcal W_k$.  Its first sum becomes a sum of squared
absolute values on $(v,\bar v)$ and is therefore semipositive; discarding it proves the final inequality.
\end{proof}

\Needspace{10\baselineskip}
\subsection{Partial-sum estimates}

\begin{prop}

\label{thm:traced-sum}
Let $K_1,K_2>0$. Suppose that
\[
 \operatorname{tr}_{\omega_g}R^{E_1}\ge K_1\varphi^*H
\]
and, for every $y\in N$, $w\in T_y^{1,0}N$, and $\eta\in(E_2)_y$,
\[
 R^{E_2}(w,\bar w,\eta,\bar\eta)
 \le K_2|w|_h^2|\eta|_H^2.
\]
Then, for $1\le k\le q$,
\[
 \sup_M\sigma_k\le\frac{k'K_2E_f}{K_1},
 \qquad k'=\min\{k,r_\varphi\}.
\]
\end{prop}

\begin{proof}
Set $s=\max_M\sigma_k$. If $s=0$, there is nothing to prove.
Otherwise, choose a maximum point $x_0$ and singular frames there.
The local support satisfies
\[
 \mathcal U_k(x)\le\sigma_k(x)\le s=\mathcal U_k(x_0).
\]
Thus $\mathcal U_k$ has a local maximum at $x_0$.
At most $k'$ of $\mu_1,\ldots,\mu_k$ are nonzero, so $k'>0$ and
$s^2\le k'\sum_{a=1}^k \mu_a^2$ by Cauchy--Schwarz.
Tracing \eqref{Ukk} and applying the curvature bounds gives
\[
 0\ge\Delta_g\mathcal U_k(x_0)
 \ge K_1\sum_{a=1}^k \mu_a^2-K_2E_fs
 \ge\frac{K_1}{k'}s^2-K_2E_fs.
\]
Dividing by $s$ proves the estimate.
\end{proof}

\begin{remark}
The estimate for $k=1$ also implies the general case, since
$\sigma_k\le k'\sigma_1$. The same argument gives two variants.

First, suppose that at every $x\in M$ there is a nonzero
$v\in T_x^{1,0}M$ such that, for all $\xi\in(E_1)_x$ and
$\eta\in(E_2)_{f(x)}$,
\[
 \begin{aligned}
 R^{E_1}(v,\bar v,\xi,\bar\xi)
 &\ge K_1|v|_g^2|\varphi(\xi)|_H^2,\\
 R^{E_2}(f_*v,\overline{f_*v},\eta,\bar\eta)
 &\le K_2|v|_g^2|\eta|_H^2.
 \end{aligned}
\]
Then $\sup_M\sigma_k\le k'K_2/K_1$, by evaluating \eqref{Ukk} on this direction at a maximum point of $\sigma_k$.

Second, retain $\operatorname{tr}_{\omega_g}R^{E_1}\ge K_1\varphi^*H$ and assume
that, for every $x\in M$, $v\in T_x^{1,0}M$, and
$\eta\in(E_2)_{f(x)}$,
\[
 R^{E_2}(f_*v,\overline{f_*v},\eta,\bar\eta)
 \le K_2|v|_g^2|\eta|_H^2.
\]
Then $\sup_M\sigma_k\le k'r_fK_2/K_1$. At a maximum point $x_0$, take a $g$-unitary basis whose last $m-t$
vectors span $\ker df_{x_0}$, where $t=\operatorname{rank}_{\mathbb C}df_{x_0}\le r_f$.
Their target-curvature terms vanish, so tracing \eqref{Ukk} gives the stated factor $r_f$.
\end{remark}

\subsection{Partial products and exterior powers}

In the formula for $\log\mathcal W_k$, the curvature terms appear without singular-value weights. This allows the estimates to use
curvature traces over $k$-dimensional fiber subspaces in place of bounds on individual fiber directions.

\begin{theorem}\label{thm:product}
Fix $1\le k\le q$ and $K_1,K_2>0$. Assume that, at every $x\in M$, each $G$-orthonormal $k$-frame $e_1,\ldots,e_k$ in $(E_1)_x$ satisfies
\begin{equation}\label{eq:product-source-trace}
 \sum_{a=1}^k
 \bigl(\operatorname{tr}_{\omega_g}R^{E_1}\bigr)(e_a,\bar e_a)
 \ge K_1\sum_{a=1}^k|\varphi(e_a)|_H^2.
\end{equation}
Then the following statements hold.
\begin{enumerate}
\item Suppose that, for every $y\in N$, $w\in T_y^{1,0}N$, and $H$-orthonormal $k$-frame $s_1,\ldots,s_k$ in $(E_2)_y$,
\begin{equation}\label{eq:product-target-trace}
 \sum_{a=1}^kR^{E_2}(w,\bar w,s_a,\bar s_a)
 \le kK_2|w|_h^2.
\end{equation}
Then
\[
 \sup_M\pi_k\le\left(\frac{K_2E_f}{K_1}\right)^k.
\]
\item Suppose that, for every $x\in M$, $v\in T_x^{1,0}M$, and $H$-orthonormal $k$-frame $s_1,\ldots,s_k$ in $(E_2)_{f(x)}$,
\begin{equation}\label{eq:product-pulled-target-trace}
 \sum_{a=1}^kR^{E_2}(f_*v,\overline{f_*v},s_a,\bar s_a)
 \le kK_2|v|_g^2.
\end{equation}
Then
\[
 \sup_M\pi_k\le\left(\frac{r_fK_2}{K_1}\right)^k.
\]
\end{enumerate}
\end{theorem}

\begin{remark}\label{rem:product-pointwise}
Each fiber sum in the theorem is independent of the orthonormal basis
of its $k$-plane. In particular,
$\operatorname{tr}_{\omega_g}R^{E_1}\ge K_1\varphi^*H$ implies
\eqref{eq:product-source-trace}. The respective bounds
\begin{align*}
 R^{E_2}(w,\bar w,\eta,\bar\eta)
 &\le K_2|w|_h^2|\eta|_H^2,\\
 R^{E_2}(f_*v,\overline{f_*v},\eta,\bar\eta)
 &\le K_2|v|_g^2|\eta|_H^2
\end{align*}
imply \eqref{eq:product-target-trace} and
\eqref{eq:product-pulled-target-trace}. Thus both estimates also hold
under these pointwise assumptions. For $k>1$, the fiber-trace conditions
allow negative values of $\operatorname{tr}_{\omega_g}R^{E_1}$
on individual fiber directions.
\end{remark}

\begin{proof}[Proof of Theorem~\ref{thm:product}]
If $r_\varphi<k$, then $\pi_k\equiv0$.  We may assume $r_\varphi\ge k$. Let $x_0$ be a global maximum point of $\pi_k$;
the maximum is positive. The local support satisfies
\[
 \mathcal W_k(x)\le\pi_k(x)\le\pi_k(x_0)=\mathcal W_k(x_0).
\]
After shrinking the neighborhood, $\mathcal W_k>0$ there, so $\log\mathcal W_k$ has a local maximum at $x_0$.
Choose singular frames $e_a,s_a$ at $x_0$ as in Proposition~\ref{prop:local-bochner}, and let
$\{\partial_1,\ldots,\partial_m\}$ be a $g$-unitary basis.
Tracing the logarithmic support formula and discarding its nonnegative
derivative term gives
\begin{align*}
 0&\ge\Delta_g\log\mathcal W_k(x_0)\\
 &\ge\sum_{i=1}^{m}\sum_{a=1}^{k}
 \left\{R^{E_1}(\partial_i,\bar\partial_i,e_a,\bar e_a)
 -R^{E_2}(f_*\partial_i,\overline{f_*\partial_i},s_a,\bar s_a)\right\}.
\end{align*}
Applying \eqref{eq:product-source-trace} and
\eqref{eq:product-target-trace} to these frames yields
\[
 K_1\sum_{a=1}^{k}\mu_a
 \le kK_2\operatorname{tr}_{\omega_g}f^*\omega_h(x_0).
\]
The arithmetic--geometric mean inequality now gives
\[
 \pi_k(x_0)^{1/k}
 \le\frac1k\sum_{a=1}^k\mu_a
 \le\frac{K_2}{K_1}
 \operatorname{tr}_{\omega_g}f^*\omega_h(x_0)
 \le\frac{K_2E_f}{K_1}.
\]
This proves part~(1).
For part~(2), put $t=\operatorname{rank}_{\mathbb C}df_{x_0}$ and choose the last $m-t$ basis vectors in $\ker df_{x_0}$, whose complex dimension is $m-t$. Their target-curvature terms vanish. Summing \eqref{eq:product-pulled-target-trace} over the remaining $t$ unit vectors therefore gives
\[
 K_1\sum_{a=1}^k\mu_a\le ktK_2\le kr_fK_2.
\]
The same arithmetic--geometric mean inequality yields
\[
 \pi_k(x_0)^{1/k}\le\frac{r_fK_2}{K_1}.
\]
\end{proof}

Since $\pi_k=\|\wedge^k\varphi\|_{\mathrm{op}}^2$ for the induced Hermitian metrics, Theorem~\ref{thm:product} also gives metric bounds on the exterior powers.

\begin{example}
\label{ex:product-mixed-sign}
Let $M=N=\mathbb{CP}^1$, $g=h$, and $f=\operatorname{id}$.
Normalize the Fubini--Study metric and the induced metrics $G_d$ on $\mathcal O(d)$ so that
\[
 R^{(\mathcal O(d),G_d)}(v,\bar v,\xi,\bar\xi)
 =d|v|_g^2|\xi|_{G_d}^2.
\]
On $E=\mathcal O(-1)\oplus\mathcal O(3)$, set
\[
 G=G_{-1}\oplus G_3,
 \qquad H=\tfrac32G_{-1}\oplus\tfrac12G_3,
 \qquad\varphi=\operatorname{id}.
\]
For $k=2$, both sides of \eqref{eq:product-source-trace} equal $2$ when $K_1=1$. The target fiber trace is $2|w|_g^2$, since constant
rescaling of each summand preserves its curvature endomorphism.
Hence both target conditions hold with $K_2=1$, and $E_f=r_f=1$.
The squared singular values are $3/2$ and $1/2$, so
\[
 \pi_2=\tfrac34\le1,\qquad \mu_1=\tfrac32>1.
\]
The fiber-trace assumptions therefore give the product estimate although the operator estimate with the same constant fails.
\end{example}

For a holomorphic map $f:(M^m,g)\to(N^n,h)$ with $m\le n$, take
$E_1=T^{1,0}M$, $E_2=T^{1,0}N$, and $\varphi=df$. The top product
is the squared complex $m$-dilation:
\[
 f^*(\omega_h^m)=\pi_m\omega_g^m.
\]
Thus any bound $\pi_m\le C^m$ gives a volume-form comparison.
If $M,N$ are compact and $m=n$, integration also gives
\[
 \deg(f)\operatorname{Vol}(N,h)
 \le C^m\operatorname{Vol}(M,g),
 \qquad
 \operatorname{Vol}(M,g)=\frac1{m!}\int_M\omega_g^m.
\]

\begin{remark}[Sharpness]\label{rem:product-sharp}
Both estimates in Theorem~\ref{thm:product} are sharp.
Take the direct-sum construction in Remark~\ref{rem:sharp-complete}
with $r=k$. Then
\[
 E_f=r_f=1,\qquad K_1=\kappa/c,\qquad K_2=\kappa,
 \qquad \pi_k=c^k.
\]
Thus equality holds in both estimates.
\end{remark}

\section{Rigidity and rank bounds}\label{sec:rank}

The rank criteria use curvature traces on complex fiber subspaces.
Let $E$ have rank $r$ over a Hermitian manifold of complex dimension $m$, and fix $1\le l\le r$ and $1\le k\le m$.
If $\sigma\subset E_x$ is an $l$-plane with a unitary basis $\{\xi_\alpha\}_{\alpha=1}^{l}$, define
\[
 R^E(v;\sigma)=\sum_{\alpha=1}^{l}
 R^E(v,\bar v,\xi_\alpha,\bar\xi_\alpha).
\]
If $\Sigma\subset T_x^{1,0}M$ is a $k$-plane with a unitary basis
$\{v_i\}_{i=1}^{k}$, put
\[
 R^E(\Sigma;\sigma)=\sum_{i=1}^{k}R^E(v_i;\sigma).
\]

\subsection{Parallelism and determinant rigidity}

We first record a parallelism statement under curvature sign conditions in every base and fiber direction. Related equality cases for Schwarz lemmas are studied in
\cite{CN22}.

\begin{prop}\label{thm:parallel}\label{thm:parallel-main}
Let $M$ be compact and connected. Suppose that for every
$x\in M$, $v\in T_x^{1,0}M$, $\xi\in(E_1)_x$, and
$\eta\in(E_2)_{f(x)}$,
\[
 R^{E_1}(v,\bar v,\xi,\bar\xi)\ge0,
 \qquad
 R^{E_2}(f_*v,\overline{f_*v},\eta,\bar\eta)\le0.
\]
Then $\nabla\varphi=0$. Its singular values and rank are constant.
\end{prop}

\begin{proof}
Consider the global smooth function
$U_1=\operatorname{tr}_G\varphi^*H$.  The $p=1$ Bochner formula gives,
at a singular-value frame,
\[
 \partial\bar\partial U_1(v,\bar v)
 =|\nabla_v\varphi|^2+
 \sum_{a=1}^{q}\mu_a
 \left\{R^{E_1}(v,\bar v,e_a,\bar e_a)
 -R^{E_2}(f_*v,\overline{f_*v},s_a,\bar s_a)\right\}.
\]
The right-hand side is nonnegative, so $U_1$ is plurisubharmonic. A plurisubharmonic function on a compact connected complex manifold is
constant.  The displayed identity therefore forces $|\nabla_v\varphi|=0$ for every $v$, and hence $\nabla\varphi=0$.
Parallel transport is unitary for the Chern connections, so it preserves
all singular values.
\end{proof}

For the top exterior power, curvature traces on fiber subspaces
suffice. The following consequence of Proposition~\ref{thm:parallel}
also admits a proof by Ni's logarithmic volume method
\cite[Corollary~3.2]{Ni21}.

\begin{theorem}\label{thm:det-rigidity}
Let $M$ be compact and connected, and assume $r_1\le r_2$.
Suppose that, for every $x\in M$, $v\in T_x^{1,0}M$, and complex
$r_1$-plane $\sigma'\subset(E_2)_{f(x)}$,
\[
 R^{E_1}(v;(E_1)_x)\ge0,
 \qquad R^{E_2}(df_x(v);\sigma')\le0.
\]
Then $\wedge^{r_1}\varphi$ is parallel. Consequently, $\pi_{r_1}=\|\wedge^{r_1}\varphi\|_{\mathrm{op}}^2$ is constant,
and full rank at one point implies full rank everywhere.
\end{theorem}

\begin{proof}
Set
\[
 L=\det E_1,\qquad F=f^*(\wedge^{r_1}E_2),
\]
with their induced metrics. We verify the curvature inequalities
in Proposition~\ref{thm:parallel} for
$\psi=\wedge^{r_1}\varphi:L\to F$.

Fix $x\in M$ and $v\in T_x^{1,0}M$.
The curvature of the determinant metric is the fiber trace of
the original curvature. Thus, for every $\xi\in L_x$,
\[
 R^L(v,\bar v,\xi,\bar\xi)
 =R^{E_1}(v;(E_1)_x)|\xi|^2\ge0.
\]

Choose an $H$-unitary basis $s_1,\ldots,s_{r_2}$ diagonalizing
$R^{E_2}(df_x(v),\overline{df_x(v)},\cdot,\cdot)$,
with eigenvalues $\theta_1,\ldots,\theta_{r_2}$.
For each increasing multi-index $I=(i_1,\ldots,i_{r_1})$, put
\[
 s_I=s_{i_1}\wedge\cdots\wedge s_{i_{r_1}},
 \qquad
 \sigma_I=\operatorname{span}\{s_{i_1},\ldots,s_{i_{r_1}}\}.
\]
The induced curvature on $F_x$ is diagonal in the basis
$\{s_I\}$, with eigenvalue
$\theta_{i_1}+\cdots+\theta_{i_{r_1}}$ on $s_I$.
Consequently, every $\eta=\sum_I\eta_I s_I\in F_x$ satisfies
\[
 \begin{aligned}
 R^F(v,\bar v,\eta,\bar\eta)
 &=\sum_I
   (\theta_{i_1}+\cdots+\theta_{i_{r_1}})|\eta_I|^2\\
 &=\sum_I R^{E_2}(df_x(v);\sigma_I)|\eta_I|^2
 \le0.
 \end{aligned}
\]
Proposition~\ref{thm:parallel} now implies that $\psi$ is parallel,
so its operator norm is constant.
Since $\psi_x\ne0$ precisely when $\varphi_x$ has rank $r_1$,
full rank at one point implies full rank everywhere.
\end{proof}

\subsection{Rank thresholds}

\begin{proof}[Proof of Theorem~\ref{thm:intro-rank}]
If $r_\varphi\ge l$, the continuous function $\pi_l$ has a positive
maximum at some $x_0\in M$. Choose unitary singular frames there.
Its support $\mathcal W_l$ is positive at $x_0$ and, after shrinking
the neighborhood, remains positive. Hence $\log\mathcal W_l$ has a
local maximum at $x_0$. Set
\[
 \sigma=\operatorname{span}\{e_1,\ldots,e_l\},\qquad
 \sigma'=\operatorname{span}\{s_1,\ldots,s_l\}.
\]
The two sign conditions, with at least one strict inequality, give a
nonzero $v\in T_{x_0}^{1,0}M$ with
$R^{E_1}(v;\sigma)-R^{E_2}(f_*v;\sigma')>0$.
Proposition~\ref{prop:local-bochner} gives
\[
 0\ge\partial\bar\partial\log \mathcal W_l(v,\bar v)
 \ge R^{E_1}(v;\sigma)-R^{E_2}(f_*v;\sigma')>0,
\]
a contradiction.
\end{proof}

The same proof applies if, for every $x$ and every pair
$(\sigma,\sigma')$, some nonzero $v\in T_x^{1,0}M$ makes
$R^{E_1}(v;\sigma)-R^{E_2}(df_x(v);\sigma')$ strictly positive.
The separate sign conditions are sufficient for this.

In particular, if $M$ is compact, then $\varphi\equiv0$ provided that,
at each $x$, there is a nonzero
$v\in T_x^{1,0}M$ for which either
$R^{E_1}(v,\bar v,\cdot,\cdot)$ is positive definite and
$R^{E_2}(df_x(v),\overline{df_x(v)},\cdot,\cdot)$ is negative
semidefinite, or the first form is positive semidefinite and the
second is negative definite. This follows from the case $l=1$ of
Theorem~\ref{thm:intro-rank}.

\subsection{Dimension and Griffiths curvature criteria}

Following Ping Li \cite[Definition~6.3]{LiP22}, we say
that $E$ is \emph{Griffiths $(k,l)$-positive} if
$R^E(\Sigma;\sigma)>0$ at every point and for every pair
$(\Sigma,\sigma)$; the other signs are defined analogously. For $k=l=1$, this is Griffiths
positivity \cite{Gr69}.

For the dimension criteria, we need a sign condition on a whole base
plane. We call $E$ \emph{strongly $RC(k,l)$-positive} if, at every
$x$ and for every fiber $l$-plane $\sigma$, there is a base $k$-plane
$\Sigma$ such that
\[
 R^E(v;\sigma)>0\qquad\text{for every }0\ne v\in\Sigma.
\]
The plane $\Sigma$ may depend on $\sigma$. The other signs have the
same quantifiers. This condition implies Li's $RC(k,l)$-positivity
\cite[Definition~6.5(1)]{LiP22}, which requires only
$R^E(\Sigma;\sigma)>0$. It differs from his \emph{uniform}
$RC(k,l)$-positivity \cite[Definition~6.5(2)]{LiP22}: there one plane
$\Sigma$ works for every $\sigma$, but positivity is required only
after taking the trace over $\Sigma$. For $l=1$, our strong condition
is the $k$-RC positivity of \cite{XYY}.
For background on RC-positivity and its geometric applications, see
\cite{Yang18,Yang21,Yang24}.

The next argument uses the kernel dimension of $df_x$, as in
\cite[Theorem~1.20]{XYY}, with fiberwise traces over $l$-planes. Put
\[
 r_{\min}(f)=\min_{x\in M}\operatorname{rank}_{\mathbb C}df_x.
\]

\Needspace{13\baselineskip}
\begin{corollary}\label{cor:dimension}
Let $M$ be compact.  Assume that $1\le l\le q$,
$1\le k\le\dim_{\mathbb C}M$, and
$1\le s\le\dim_{\mathbb C}N$.
\begin{enumerate}
\item If $E_1$ is strongly $RC(k,l)$-positive, $E_2$ is strongly
$RC(s,l)$-nonpositive, and $k+s>\dim_{\mathbb C}N$, then
$r_\varphi<l$.
\item If $E_1$ is strongly $RC(k,l)$-nonnegative, $E_2$ is strongly
$RC(s,l)$-negative, and
\[
 k+s>\dim_{\mathbb C}M+\dim_{\mathbb C}N-r_{\min}(f),
\]
then $r_\varphi<l$.
\end{enumerate}
\end{corollary}

\begin{proof}
Fix $x\in M$ and complex $l$-planes
\[
 \sigma\subset(E_1)_x,
 \qquad \sigma'\subset(E_2)_{f(x)}.
\]
Set
\[
 V=T_x^{1,0}M,\qquad V'=T_{f(x)}^{1,0}N,\qquad
 L=df_x:V\longrightarrow V'.
\]
Write
\[
 m=\dim_{\mathbb C}M,\qquad n=\dim_{\mathbb C}N,\qquad
 r_x=\operatorname{rank}_{\mathbb C}L,\qquad K_x=\ker L.
\]
Then $\dim_{\mathbb C}K_x=m-r_x$.

For part~(1), choose complex subspaces
$\Sigma\subset V$ and $\Pi\subset V'$ supplied by the two strong $RC$
assumptions, with
$\dim_{\mathbb C}\Sigma=k$ and $\dim_{\mathbb C}\Pi=s$.
If $\Sigma\cap K_x\ne\{0\}$, choose
$0\ne v\in\Sigma\cap K_x$.  Then
\[
 R^{E_1}(v;\sigma)>0,\qquad
 R^{E_2}(df_x(v);\sigma')=R^{E_2}(0;\sigma')=0.
\]
Theorem~\ref{thm:intro-rank} applies.

Suppose that $\Sigma\cap K_x=\{0\}$.  Then $L|_\Sigma$ is injective, so
$\dim_{\mathbb C}L(\Sigma)=k$.  The Grassmann dimension inequality gives
\[
 \dim_{\mathbb C}\bigl(L(\Sigma)\cap\Pi\bigr)
 \ge k+s-n>0.
\]
Choose $0\ne w\in L(\Sigma)\cap\Pi$ and $v\in\Sigma$ with $L(v)=w$.
The two curvature traces again have the signs required by
Theorem~\ref{thm:intro-rank}.

For part~(2), choose $\Sigma$ and $\Pi$ from the strong nonnegative and
strong negative assumptions.  Since
\[
 \dim_{\mathbb C}L(\Sigma)
 =k-\dim_{\mathbb C}(\Sigma\cap K_x)
 \ge k-m+r_x,
\]
we have
\begin{align*}
 \dim_{\mathbb C}\bigl(L(\Sigma)\cap\Pi\bigr)
 &\ge \dim_{\mathbb C}L(\Sigma)+s-n\\
 &\ge k+s-m-n+r_x\\
 &\ge k+s-m-n+r_{\min}(f)>0.
\end{align*}
Choose $0\ne w\in L(\Sigma)\cap\Pi$ and $v\in\Sigma$ with $L(v)=w$.
Then
\[
 R^{E_1}(v;\sigma)\ge0,
 \qquad R^{E_2}(df_x(v);\sigma')<0.
\]
Apply Theorem~\ref{thm:intro-rank}.
\end{proof}

A trace over base directions can also supply the common direction
required by Theorem~\ref{thm:intro-rank}.

\begin{corollary}\label{thm:averaged}
Let $M$ be compact, set $m=\dim_{\mathbb C}M$ and
$n=\dim_{\mathbb C}N$, and let $1\le l\le q$.
\begin{enumerate}
\item If $E_1$ is Griffiths $(m,l)$-positive and $E_2$ is
Griffiths $(1,l)$-nonpositive, then $r_\varphi<l$.
\item Suppose that $E_1$ is Griffiths $(1,l)$-positive.
If $m>n$, then $r_\varphi<l$.
If $m\le n$, the same conclusion holds provided that $E_2$
is Griffiths $(m,l)$-nonpositive.
\end{enumerate}
\end{corollary}

\begin{proof}
Fix $x\in M$ and complex $l$-planes
$\sigma\subset(E_1)_x$ and $\sigma'\subset(E_2)_{f(x)}$.
For part~(1), the positive trace
$R^{E_1}(T_x^{1,0}M;\sigma)$ gives a vector $v$ in a unitary basis
with $R^{E_1}(v;\sigma)>0$. The target trace is nonpositive on
$df_x(v)$, so Theorem~\ref{thm:intro-rank} applies.

For part~(2), if $\ker df_x\ne\{0\}$, choose a nonzero vector
$v$ in this kernel. Its curvature trace in $E_1$ is positive, while the
pulled-back target trace is zero. Otherwise $df_x$ is injective, so $m\le n$.
Choose an $h$-unitary basis $w_1,\ldots,w_m$ of
$\operatorname{im}df_x$. The target hypothesis gives
\[
 \sum_{j=1}^{m}R^{E_2}(w_j;\sigma')\le0.
\]
Thus some $w_j$ has nonpositive target trace. Its nonzero preimage
under $df_x$ has positive curvature trace in $E_1$. In either case,
Theorem~\ref{thm:intro-rank} gives $r_\varphi<l$.
\end{proof}

\begin{example}[Sharpness with mixed-sign curvature]\label{ex:rank-two}
Let $M=N=\mathbb{CP}^1$ with the Fubini--Study metric and let
$f=\operatorname{id}_{\mathbb{CP}^1}$. Equip
$L=\mathcal O_{\mathbb{CP}^1}(1)$ with its Fubini--Study metric and set
\[
 E_1=L^{-1}\oplus L^{\otimes2}\oplus L^{\otimes2},
 \qquad E_2=\mathcal O^{\oplus2},
\]
with the direct-sum metric on $E_1$ and the flat metric on $E_2$.

For a unit vector $\tau\in T_x^{1,0}\mathbb{CP}^1$, the Hermitian
endomorphism associated with
$R^{E_1}(\tau,\bar\tau,\cdot,\cdot)$ has eigenvalues
\[
 -\kappa,\qquad 2\kappa,\qquad 2\kappa
\]
for a constant $\kappa>0$.  Thus its trace on every complex two-plane
$\sigma\subset(E_1)_x$ is at least
$-\kappa+2\kappa=\kappa>0$.  The target is flat.  Therefore
Theorem~\ref{thm:intro-rank}, with $l=2$, gives
$\operatorname{rank}_{\mathbb C}\varphi<2$ for every holomorphic bundle
map $\varphi:E_1\to E_2$ covering $f$.

The bound is sharp.  Let $s_0=Z_0$ and $s_1=Z_1$ be the standard
homogeneous-coordinate sections of $L$.  They have no common zero.  For
$u\in L_x^{-1}$ and $v_1,v_2\in L_x^{\otimes2}$, define
\[
 \varphi_x(u,v_1,v_2)
 =\bigl(u(s_0(x)),u(s_1(x))\bigr)\in\mathbb C^2.
\]
Since $s_0$ and $s_1$ have no common zero,
$\operatorname{rank}_{\mathbb C}\varphi_x=1$ for every $x$.

The negative $L^{-1}$ direction prevents linewise positivity, while
the two-plane trace is positive. The exterior powers satisfy
\[
 \wedge^2E_1=L\oplus L\oplus L^{\otimes4},
 \qquad \wedge^2E_2=\mathcal O.
\]
Thus the same rank bound also follows from
\cite[Theorem~1.17]{XYY} applied to $\wedge^2\varphi$.
\end{example}

\section*{Acknowledgements}
The authors thank Professor L. Ni for valuable comments on
Schwarz estimates and related suggestions, and Professors
S. K. Wang, K. Wu, X. A. Zheng, X. Y. Zhou, and H. W. Xu for their
encouragement. This work was supported by the National Natural Science Foundation of China (Grant No.~11901164) and the Young Backbone Teacher Fund of Henan Polytechnic University (Grant No.~2025XQG-12).

\bibliographystyle{amsplain}

\bigskip
\begingroup
\small
\raggedright
\setlength{\parindent}{0pt}
\setlength{\parskip}{0pt}

\Needspace{5\baselineskip}
\textsc{Xu Pang, School of Mathematical Sciences,
Guizhou Normal University, Guiyang 550025, China}\par
\nobreak
\textit{Email address:} \texttt{2531528680@qq.com}\par

\medskip
\Needspace{5\baselineskip}
\textsc{Yao Sha, School of Mathematics and Information Science,
Henan Polytechnic University, Jiaozuo 454003, China}\par
\nobreak
\textit{Email address:} \texttt{yaosha@hpu.edu.cn}\par

\medskip
\Needspace{5\baselineskip}
\textsc{Ren Xin-An, Department of Mathematics,
China University of Mining and Technology,
Xuzhou 221116, China}\par
\nobreak
\textit{Email address:} \texttt{renx@cumt.edu.cn}\par

\endgroup

\end{document}